\title{An index inequality for embedded pseudoholomorphic curves in
 symplectizations}
\author{Michael Hutchings}
\date{}

\documentclass[12pt]{article}

\usepackage{amssymb}
\usepackage{latexsym}
\usepackage{amsmath}
\usepackage{theorem}

\usepackage{mathrsfs}
\newcommand{\mc}[1]{{\mathscr #1}}

\newtheorem{theorem}{Theorem}[section]
\newtheorem{proposition}[theorem]{Proposition}

\newtheorem{lemma}[theorem]{Lemma}
\newtheorem{assumption}[theorem]{Assumption}

{\theorembodyfont{\rmfamily}
\newtheorem{definition}[theorem]{Definition}
\newtheorem{remark}[theorem]{Remark}

\newtheorem{example}[theorem]{Example}
}

\newenvironment{proof}
{\noindent {\em Proof.}}
{\hfill $\Box$}

\newcommand{\qed} {\hfill$\Box$}

\newcommand{\pin}{p_{\op{in}}}
\newcommand{\pout}{p_{\op{out}}}
\newcommand{\murel}{I}

\newcommand{\C}{{\mathbb C}}
\newcommand{\Q}{{\mathbb Q}}
\newcommand{\R}{{\mathbb R}}
\newcommand{\Z}{{\mathbb Z}}
\newcommand{\op}{\operatorname}
\newcommand{\dbar}{\overline{\partial}}

\newcommand{\End}{\op{End}}

\newcommand{\Hom}{\op{Hom}}

\newcommand{\tensor}{\otimes}

\newcommand{\union}{\bigcup}

\newcommand{\define}[1]{{\em #1\/}}

\newcommand{\rb}[1]{\raisebox{1.5ex}[-1.5ex]{#1}}

\begin{document}

\maketitle

\begin{abstract}
Let $\Sigma$ be a surface with a symplectic form, let $\phi$ be a
symplectomorphism of $\Sigma$, and let $Y$ be the mapping torus of
$\phi$.  We show that the dimensions of moduli spaces of embedded
pseudoholomorphic curves in $\R\times Y$, with cylindrical ends
asymptotic to periodic orbits of $\phi$ or multiple covers thereof,
are bounded from above by an additive relative index.  We deduce some
compactness results for these moduli spaces.

This paper establishes some of the foundations for a program with
Michael Thaddeus, to understand the Seiberg-Witten Floer homology of
$Y$ in terms of such pseudoholomorphic curves.  Analogues of our
results should also hold in three dimensional contact topology.
\end{abstract}

\section{Introduction}

\subsection{Motivation}

Let $\Sigma$ be a compact connected surface, possibly with boundary,
with a symplectic form $\omega$.  Let $\phi$ be a
symplectomorphism of $\Sigma$, and let
\begin{equation}
\label{eqn:mappingTorus}
Y = \frac{\Sigma\times\R}{(x,t+1)\sim(\phi(x),t)}
\end{equation}
be the mapping torus of $\phi$.  In this paper we will study the
dimensions of moduli spaces of embedded pseudoholomorphic curves in
$\R\times Y$, for a suitable almost complex structure.  We have two
basic motivations for this study.

Our first motivation, in the case $\partial\Sigma=\emptyset$, is to
understand the Seiberg-Witten Floer homology of the three-manifold
$Y$, an invariant which counts solutions to the Seiberg-Witten
equations on the four-manifold $\R\times Y$.  
Taubes has shown that counting solutions to the Seiberg-Witten
equations on a closed symplectic four-manifold is equivalent to
appropriately counting embedded pseudoholomorphic curves in it
\cite{t1,t2}.  The noncompact 4-manifold $\R\times Y$ has a natural
symplectic form, and it is plausible that a version of Taubes's
theorem should hold here, relating the Seiberg-Witten Floer homology
of $Y$ to a new version of symplectic Floer homology, which we call
``periodic Floer homology''.  This is the homology of a chain complex
in which the chains are generated by unions of periodic orbits of
$\phi$ with multiplicity, where hyperbolic orbits cannot have
multiplicity greater than one.  The differential counts certain
pseudoholomorphic curves in $\R\times Y$, which we call ``flow
lines''; see \S\ref{sec:flowLines} for the precise definition.  Flow
lines need not have genus zero; but they must be embedded, except that
there may be repeated ``trivial cylinders'' which do not intersect any
other component.  (The precise definition of periodic Floer homology
and the conjectured relation with Seiberg-Witten Floer homology are
explained in \cite{ht1}.  This relation is also suggested by the
approach of Salamon \cite{sa} involving symmetric products of
$\Sigma$.  While Seiberg-Witten Floer homology is defined for closed
oriented 3-manifolds, one might also be able to define a version of it
for a compact oriented 3-manifold whose boundary is a union of tori,
e.g.\ when $\partial\Sigma\neq\emptyset$, cf.\ \cite{t3}.)

Some analytic work is needed to show that periodic Floer homology is
well defined.  A first step is to compute the dimensions of the moduli
spaces of flow lines.  This index calculation, which turns out to be
rather involved, is the main content of the present paper.  In
addition to determining the grading of periodic Floer homology, the
index theory leads to compactness results for these moduli spaces.
The compactness results allow one to define the differential $\delta$
in periodic Floer homology as a certain count of flow lines, and are a
step towards proving that $\delta^2=0$ and that the homology is
invariant under appropriate isotopy of $\phi$.

A second motivation for this paper is provided by three-dimensional
contact topology.  It is interesting to consider embedded
pseudoholomorphic curves in the symplectization of a contact
3-manifold, as such curves have found important topological
applications in the work of Hofer, Wysocki, and Zehnder.  We expect
that analogues of our index theorem and compactness will hold in this
setting, see \S\ref{sec:conclusion}.  One should further be able to
define an analogue of periodic Floer homology on a contact
three-manifold.  This would be a possibly interesting variant of the
symplectic field theory of Eliashberg, Givental, and Hofer \cite{sft}.
The main difference between periodic Floer homology and symplectic
field theory is that the latter theory counts pseudoholomorphic curves
that are not necessarily embedded.

\subsection{Flow lines}
\label{sec:flowLines}

We now give the precise definition of the pseudoholomorphic curves
that we will study.  We begin with their boundary values.  A
\define{periodic orbit} of $\phi$ is a finite set $\gamma$ of points
in $\Sigma$ which are cyclically permuted by $\phi$.  If $x\in\gamma$
and $p=|\gamma|$ is the period of $\gamma$, the linearized
Poincar\'{e} return map of $\phi^p$ at $x$ is a symplectic linear map
$d\phi^p:T_x\Sigma\to T_x\Sigma$.  We say that $\gamma$ is
\define{nondegenerate} if $1-d\phi^{kp}$ is invertible for all
positive integers $k$.  Except where stated otherwise, we assume that
all periodic orbits of $\phi$ are nondegenerate; this holds for
generic $\phi$.  We define the
\define{Lefschetz sign}
\[
(-1)^{\epsilon(\gamma)}= \op{sign}\det(1-d\phi^p)\in\{+1,-1\}.
\]
We say that the orbit $\gamma$ is \define{hyperbolic} if $d\phi^p$ has
real eigenvalues, and \define{elliptic} otherwise.  The Lefschetz sign
is $-1$ if $\gamma$ is hyperbolic with positive eigenvalues (of
$d\phi^p$); and the Lefschetz sign is $+1$ if $\gamma$ is elliptic or
hyperbolic with negative eigenvalues.

A periodic orbit $\gamma$ determines an embedded oriented circle in
$Y$.  More precisely, by equation \eqref{eqn:mappingTorus}, the
mapping torus $Y$ fibers over $S^1=\R/\Z$.  The flow in the $\R$
direction of $\Sigma\times\R$ induces a vector field on $Y$, which we
denote by $\partial_t$, and which can also be regarded as a connection
on the bundle $Y\to S^1$.  Periodic orbits as defined above correspond
to embedded closed orbits of the vector field $\partial_t$.  From now
on we will identify periodic orbits with the corresponding circles in
$Y$.

\begin{definition}
An \define{orbit set} is a finite set of pairs
$\alpha=\{(\alpha_1,m_1),\ldots,(\alpha_k,m_k)\}$, where
$\alpha_1,\ldots,\alpha_k$ are disjoint periodic orbits, and
$m_1,\ldots,m_k$ are positive integers (``multiplicities'').  The
orbit set $\alpha$ is \define{admissible} if $m_i=1$ whenever
$\alpha_i$ is hyperbolic.
\end{definition}
(Admissible orbit sets generate the chains
in periodic Floer homology, and the mod 2 grading is given by the
product of the Lefschetz signs; see \cite{ht1}.)  If $\alpha_i$ has
period $p_i$, we define the \define{degree} of $\alpha$ to be
$d=\sum_im_ip_i$.

To discuss pseudoholomorphic curves, we need to specify an almost
complex structure on $\R\times Y$.  Let $E$ denote the vertical
tangent bundle of $Y\to S^1$.  The symplectic form $\omega$ on
$\Sigma$ defines a symplectic structure on $E$, and we choose an
$\omega$-compatible almost complex structure $J$ on $E$.
Compatibility here just means that $J$ respects the orientation, i.e.\
$\omega(v,Jv)\ge 0$.  We extend $J$ to an $\R$-invariant almost
complex structure on the four-manifold $\R\times Y$, which we also
denote by $J$, by specifying that
\[
J(\partial_s)=\partial_t,
\]
where $s$ denotes the $\R$ coordinate.  We call an almost complex
structure $J$ on $\R\times Y$ obtained this way
\define{admissible}.  We assume throughout the paper that $J$ is
admissible, except in Theorem~\ref{thm:compactness}(b) and
\S\ref{sec:almostAdmissible}.

We consider pseudoholomorphic curves in $\R\times Y$ of the form
$(C,j,u)$, where $C$ is a punctured compact Riemann surface and each
component has at least two punctures; $j$ is a complex structure on
$C$; and $u:C\to\R\times Y$ is pseudoholomorphic, $u_*\circ j=J\circ
u_*$.  We mod out by reparametrization, i.e.\ we declare $(C,j,u)$ and
$(C',j',u')$ to be equivalent if there is a diffeomorphism from $C$ to
$C'$ intertwining $j$ with $j'$ and $u$ with $u'$.  When $u$ is an
embedding (or almost an embedding), we identify $(C,j,u)$ with the
image of $u$ in $\R\times Y$, which we denote simply by $C$.

If $\gamma$ is a periodic orbit, then $\R\times{\gamma}$ is a
pseudoholomorphic curve, which we call a \define{trivial cylinder}.
More generally, a pseudoholomorphic curve may have an end smoothly
asymptotic as $s\to +\infty$ to $\R\times{\gamma}^k$, where $s$
denotes the $\R$ coordinate and ${\gamma}^k$ denotes a $k$-fold
connected cover of ${\gamma}$; we call this an
\define{outgoing end} at $\gamma$ of multiplicity $k$.  We call an end
asymptotic to $\R\times{\gamma}^k$ as $s\to-\infty$ an
\define{incoming end}.

\begin{definition}
A \define{flow line} from the orbit set $\{(\alpha_i,m_i)\}$ to the
orbit set $\{(\beta_j,n_j)\}$ is a pseudoholomorphic curve
$C\subset\R\times Y$ as above, such that:
\begin{itemize}
\item
$C$ is embedded, except that there may be repeated
trivial cylinders, which do not intersect other components of $C$.
\item
$C$ has outgoing ends at $\alpha_i$ with total multiplicity $m_i$,
incoming ends at $\beta_j$ with total multiplicity $n_j$, and no other
ends.
\end{itemize}
\end{definition}
(The differential in periodic Floer homology is defined by a certain
count of flow lines, see \cite{ht1}.)

\subsection{Some assumptions on $\phi$ and $J$}
\label{sec:nice}

In this paper we will usually assume that $\phi$ and $J$ satisfy some
additional conditions.

\begin{definition}
Let $\gamma$ be a periodic orbit, of period $p$, containing
$x\in\Sigma$.  We say that $(\phi,J)$ is \define{admissible} near
$\gamma$ if:
\begin{description}
\item{(i)}
There exists a neighborhood $U$ of $x$ in $\Sigma$, and a symplectic
identification of $U$ with a subset of $\R^2$, sending $x$ to $0$, on
which $\phi^p$ is {\em linear} near $0$.
\item{(ii)}
The restriction to $E$ of the almost complex structure $J$ on some
tubular neighborhood $N$ of $\gamma$ is given by a {\em constant\/} matrix on
each fiber of the bundle $N\to\gamma$ induced by the projection $Y\to
S^1$, for a trivialization that is linear with respect to the
identication in (i).
\end{description}
\end{definition}

\begin{definition}
For a positive integer $d$, we say that $(\phi,J)$ is
\define{$d$-admissible} if $J$ is admissible and:
\begin{description}
\item{(a)}
 $(\phi,J)$ is admissible near all periodic orbits of
period $p\le d$.
\item{(b)}
Every boundary component of $\Sigma$ has an identification of a
neighborhood of it with $(-\epsilon,0]\times (\R/\Z)$ with coordinates
$(x,y)$, in which $J(\partial_x)=\partial_y$ and
$\phi(x,y)=(x,y+\theta)$, where $q\theta\notin \Z$ for all integers
$1\le q\le d$.
\end{description}
\end{definition}

Condition (a) will simplify the asymptotic analysis of the ends of
flow lines, allowing us to focus on the topological calculations.
This assumption can probably be removed, but to carry out our proofs
without it, one would need to generalize the asymptotic analysis of
\cite{hwz1,abbas}.  In any case, we can achieve this condition by
perturbing $\phi$ and $J$, so this assumption entails no loss of
generality for defining periodic Floer homology, see
\cite{ht1}.
%The reason is that periodic Floer homology is divided
%into ``sectors'' in which we only consider orbit sets of a given
% degree $d$.

Condition (b) ensures that a flow line can between orbit sets of
degree $d$ can never approach the boundary, by the maximum principle,
because on any flow line, $x$ is a harmonic function in a neighborhood
of the boundary of $\R\times Y$.  Note that if $\theta$ is a rational
number with denominator larger than $d$, then there will be degenerate
periodic orbits near the boundary of period greater than $d$, but
these will not arise in our discussion below.

\subsection{The index inequality}

We will obtain a bound on the dimensions of the moduli spaces of flow
lines in terms of a relative index, which we now define.

Let $\alpha=\{(\alpha_i,m_i)\}$ and $\beta=\{(\beta_j,n_j)\}$
be orbit sets.  If there are any flow lines from $\alpha$ to $\beta$,
then their total homology classes must be equal:
\begin{equation}
\label{eqn:totalHomologyClass}
\sum_im_i[{\alpha_i}]=\sum_jn_j[{\beta_j}]=h\in H_1(Y).
\end{equation}
Assuming \eqref{eqn:totalHomologyClass}, define $H_2(Y;\alpha,\beta)$
to be the set of relative homology classes of 2-chains $W$ in $Y$ with
$\partial W =
\sum_im_i{\alpha_i} - \sum_jn_j{\beta_j}$.  This is an affine space
modelled on $H_2(Y)$.  If $C$ is a flow line from $\alpha$ to $\beta$,
then its projection to $Y$ determines a class $[C]\in
H_2(Y;\alpha,\beta)$.

\begin{definition}
Given $Z\in H_2(Y;\alpha,\beta)$, we define the \define{relative
index}
\[
\begin{split}
\murel(\alpha,\beta;Z) =&\; c_1(E|_Z,\tau)+Q_\tau(Z,Z)\\
&+\sum_i\sum_{k=1}^{m_i}\mu_\tau({\alpha_i}^k)
-\sum_j\sum_{k=1}^{n_j}\mu_\tau({\beta_j}^k).
\end{split}
\]
\end{definition}
Here $\tau$ is a trivialization of $E$ over the ${\alpha_i}$'s and
${\beta_j}$'s; $c_1$ is the relative first Chern class; $Q_\tau$ is a
relative intersection pairing; and $\mu_\tau$ is the
Conley-Zehnder index.  These notions are explained in detail in
\S\ref{sec:relativeIndex}.  The relative index has the following basic
properties which are proved in \S\ref{sec:adjunction}.

\begin{proposition}[properties of the relative index]
\label{prop:relativeIndex}
\begin{description}
\item{ }
\item{(a)}
(Well defined)
$\murel(\alpha,\beta;Z)$ does not depend on $\tau$.
\item{(b)}
(Additivity)
$
\murel(\alpha,\beta;Z)+\murel(\beta,\gamma;W)=\murel(\alpha,\gamma;Z+W).
$
\item{(c)}
(Parity and Lefschetz signs)
If $\alpha$ and $\beta$ are admissible, then
\[
\murel(\alpha,\beta;Z) \equiv
\sum_i\epsilon(\alpha_i)-\sum_j\epsilon(\beta_j)\mod 2.
\]
\item{(d)}
(Change of homology class)
\[
\murel(\alpha,\beta;Z)-\murel(\alpha,\beta;Z') = \langle c(h), Z-Z'
\rangle
\]
where the ``index ambiguity class'' $c(h)$ is defined by
\[
c(h)= c_1(E)+2\op{PD}(h)\in H^2(Y;\Z).
\]
\end{description}
\end{proposition}

Let $C$ be a flow line from $\alpha$ to $\beta$, and let
$\mc{M}_C$ denote the component of the moduli space of flow lines
containing $C$.  The main result of this paper is the following
theorem, which is proved in
\S\ref{sec:beginComputation}--\S\ref{sec:cylinders}.

\begin{theorem}[index inequality]
\label{thm:index}
Let $\alpha$ and $\beta$ be orbit sets of degree $d$.  Assume that
$(\phi,J)$ is $d$-admissible and $J$ is generic.  If $C$ is a flow
line from $\alpha$ to $\beta$, then $\mc{M}_C$ is a manifold and
\begin{equation}
\label{eqn:indexInequality}
\dim(\mc{M}_C)\le\murel(\alpha,\beta;[C]).
\end{equation}
Equality holds only if $C$ is admissible in the sense of
Definition~\ref{def:admissible}.
\end{theorem}

The index theorem is an inequality, rather than an equality, because
the dimension of the moduli space depends on some additional discrete
choices.  Most importantly, the multiplicities of the outgoing ends at
$\alpha_i$ determine a partition of the integer $m_i$, and the
multiplicities of the incoming ends at $\beta_j$ determine a partition
of $n_j$.  The dimension depends in part on the genus of $C$, which is
determined by a relative adjunction formula; this formula involves the
writhes of braids determined by the ends of $C$, which in turn have
bounds depending on the above partitions.

For each periodic orbit $\gamma$ and each positive integer $m$, we
{\em a priori\/} define two partitions of $m$, the ``incoming'' and
``outgoing'' partitions, denoted by $\pin(\gamma,m)$ and
$\pout(\gamma,m)$, see \S\ref{sec:combinatorics}.  If $C$ contains no
trivial cylinders, then $C$ is ``admissible'', i.e.\ equality can hold
in the index theorem, when the partition of $m_i$ determined by the
outgoing ends at $\alpha_i$ agrees with $\pout(\alpha_i,m_i)$, and the
partition of $n_j$ determined by the incoming ends agrees with
$\pin(\beta_j,n_j)$.  When $C$ contains trivial cylinders, the
criterion for $C$ to be ``admissible'' is more complicated, but again
is phrased in terms of incoming and outgoing partitions.

\subsection{Compactness}
\label{sec:compactnessIntro}

Theorem~\ref{thm:index} is strong because the upper bound $I$ on the
dimensions of the moduli spaces is additive, by
Proposition~\ref{prop:relativeIndex}(b).  This additivity leads to
some compactness results which we now state.

If $\alpha$ and $\beta$ are orbit sets and $Z\in H_2(Y;\alpha,\beta)$,
let $\mc{M}(\alpha,\beta;Z)$ denote the moduli space of flow lines $C$
from $\alpha$ to $\beta$ with relative homology class $[C]=Z$.  Note
that $\R$ acts on $\mc{M}(\alpha,\beta;Z)$ by translation in the $\R$
direction of $\R\times Y$.  The natural flat connection on $\R\times
Y\to\R\times S^1$ allows us to extend the symplectic form $\omega$ on
the fibers to a canonical closed 2-form on $\R\times Y$, which we also
denote by $\omega$.  If $i\in\Z$ and $R\in\R$, define
\[
\mc{M}_i^R(\alpha,\beta) =
\bigcup_{\substack{\murel(\alpha,\beta;Z)=i\\\int_Z\omega<R}}
\mc{M}(\alpha,\beta;Z)/\R.
\]

\begin{theorem}[compactness]
\label{thm:compactness}
Let $\alpha$ and $\beta$ be orbit sets of degree $d$, and assume that
$(\phi,J)$ is $d$-admissible and $J$ is generic.
\begin{description}
\item{(a)}
Suppose that
$d>\op{genus}(\Sigma)$, or $d=1$, or
$\partial\Sigma\neq\emptyset$.
Then:
\begin{description}
\item{(i)}
$\mc{M}_1^R(\alpha,\beta)$ is finite for each $R$.
\item{(ii)}
If $\alpha,\beta$ are admissible, then $\mc{M}_2^R(\alpha,\beta)$ has
a compactification $\overline{\mc{M}_2^R}(\alpha,\beta)$ with a natural map
\begin{equation}
\label{eqn:brokenFlowLines}
\overline{\mc{M}_2^R}(\alpha,\beta) \setminus
\mc{M}_2^R(\alpha,\beta) \longrightarrow
\union_\gamma \mc{M}_1^R(\alpha,\gamma)\times \mc{M}_1^R(\gamma,\beta).
\end{equation}
\end{description}
\item{(b)}
In general, $J$ can be perturbed to an ``almost $d$-admissible''
almost complex structure on $\R\times Y$, see
Definition~\ref{def:almost}, such that (i) and (ii) above hold, and
the index theorem
\ref{thm:index} still holds.
\end{description}
\end{theorem}

Assertion (i) allows one to define the differential $\delta$ in
periodic Floer homology as a certain count of flow lines, and
assertion (ii) is a step towards proving that $\delta^2=0$ and similarly
that the homology is invariant under suitable isotopies, see
\cite{ht1}.  (To complete the proof that $\delta^2=0$, one needs a
gluing theorem to show that for a given $\gamma$, over a pair of flow
lines on the right side of \eqref{eqn:brokenFlowLines} with total
energy less than $R$, the map
\eqref{eqn:brokenFlowLines} is odd-to-one if $\gamma$ is admissible,
and even-to-one if $\gamma$ is not admissible, see
\cite{ht1}.)  

The proof of Theorem~\ref{thm:compactness}, given in
\S\ref{sec:multiplyCovered}--\S\ref{sec:compactness}, is mostly a
standard application of Gromov compactness.  The main difficulty is
that {\em a priori\/}, sequences of embedded pseudoholomorphic curves
could converge to non-embedded curves.  The hardest part of the proof
is to check that sequences of flow lines in the relevant moduli spaces
do not converge to multiply covered pseudoholomorphic curves.  We will
do this by enhancing the index inequality to show that an embedded (or
almost embedded) curve underlying such a multiply covered curve would
live in a moduli space of negative expected dimension, and hence does
not exist for generic $J$.  Another problem is that for an admissible
almost complex structure, the fibers of the projection $\R\times
Y\to\R\times S^1$ are pseudoholomorphic, and these might bubble off as
the complex structures on the pseudoholomorphic curves in a sequence
degenerate.  (A related issue arises in the work of Ionel-Parker
\cite{ip2} on Gromov-Witten invariants of symplectic sums.)  Under the
assumption in part (a), we will see that this possibility can be ruled
out.  If this assumption does not hold, in particular if $g>1$, then
the pseudoholomorphic curves corresponding to fibers live in moduli
spaces of negative expected dimension.  We will then see in the proof
of part (b) that an appropriate perturbation of the almost complex
structure will make such curves disappear, without interfering with
the rest of the proofs of the index theorem and compactness.

\subsection{Other results}

In \S\ref{sec:euler}, as a by-product of the index calculations, we
work out a formula for the Euler characteristic of flow lines.  This
is useful for computing periodic Floer homology in specific examples.
In \S\ref{sec:conclusion} we make some concluding remarks.

\section{The relative index}
\label{sec:relativeIndex}

We now give a detailed explanation of the relative index $\murel$ which
appears in the index theorem.  We will prove its basic properties
(Proposition~\ref{prop:relativeIndex}) in
\S\ref{sec:relativeIndexProperties}, after introducing the relative
adjunction formulas.

Let $\alpha=\{(\alpha_i,m_i)\}$
and $\beta=\{(\beta_j,n_j)\}$ be orbit sets with the same total
homology class \eqref{eqn:totalHomologyClass}.

\subsection{Trivializations}

If $\gamma$ is a periodic orbit, let $\mc{T}(\gamma)$ denote the set
of homotopy classes of symplectic trivializations of $E|_\gamma$.
This is an affine space over $\Z$.  We adopt the sign convention that
if $\tau_1,\tau_2:E|_\gamma\to S^1\times\R^2$ are two trivializations,
then $\tau_2-\tau_1$ is the degree of
$\tau_1\circ\tau_2^{-1}:S^1\to\op{Sp}(2,\R)\approx S^1$.

We say that a nonvanishing section of $E|_\gamma$ is
\define{$\tau$-trivial} if its winding number with respect to $\tau$
is zero.

We define $\mc{T}(\alpha,\beta)=
\prod_i\mc{T}(\alpha_i)\times\prod_j\mc{T}(\beta_j)$.  If
$\tau\in\mc{T}(\alpha,\beta)$, we denote the corresponding elements of
$\mc{T}(\alpha_i)$ and $\mc{T}(\beta_j)$ by $\tau_i^+$ and $\tau_j^-$.

\subsection{The relative first Chern class}

If $Z\in H_2(Y;\alpha,\beta)$ and $\tau\in\mc{T}(\alpha,\beta)$, we
define the \define{relative first Chern class} $c_1(E|_Z,\tau)\in\Z$
as follows: Choose a surface $S$ with boundary and a map $f:S\to Y$
representing $Z$, choose a $\tau$-trivial (nonvanishing) section
$\psi$ of $f^*E|_{\partial S}$, extend $\psi$ to a section of $f^*E$
over $S$, and define $c_1(E|_Z,\tau)$ to be the signed number of
zeroes of this extension.

This has the following elementary
properties, cf.\ \cite{ht1}.  First, $c_1(E|_Z,\tau)$ depends only on
$Z$ and $\tau$.  Second, if we change the homology class $Z$, then
\begin{equation}
\label{eqn:chernHom}
c_1(E|_Z,\tau)-c_1(E|_{Z'},\tau)=\langle c_1(E),Z-Z'\rangle
\end{equation}
where $c_1(E)\in H^2(Y;\Z)$ is the ordinary first Chern class.  Third,
under a change of trivialization, we have
\begin{equation}
\label{eqn:chernTriv}
c_1(E|_Z,\tau)-c_1(E|_Z,\tau')=
\sum_im_i({\tau_i'}^+-\tau_i^+)-\sum_jn_j({\tau_j'}^--\tau_j^-).
\end{equation}

\subsection{The Conley-Zehnder index}
\label{sec:CZ}

If $\gamma$ is a periodic orbit and $\tau\in\mc{T}(\gamma)$, we define
the \define{Conley-Zehnder index} $\mu_\tau(\gamma)\in\Z$ as follows.
Recall that there is a natural connection on $Y\to S^1$, whose
monodromy is given by $\phi$.  The linearized parallel transport
induces a connection on $E|_\gamma$.  If we traverse $\gamma$ once,
then with respect to the trivialization $\tau$, parallel transport
defines a path in $\op{Sp}(2,\R)$ from the identity to the linearized
return map $d\phi^p$, where $p$ is the period of $\gamma$.  We define
$\mu_\tau(\gamma)$ to be the Maslov index of this path of symplectic
matrices, see e.g. \cite{sz}.  Likewise we define $\mu_\tau(\gamma^k)$
to be the Maslov index of the path in $\op{Sp}(2,\R)$ obtained by
traveling $k$ times around $\gamma$.  Since we are assuming that all
periodic orbits are nondegenerate, these Maslov indices are defined.

The general theory of the Maslov index implies that if we change the
trivialization, then
\begin{equation}
\label{eqn:CZTriv}
\mu_\tau(\gamma^k)-\mu_{\tau'}(\gamma^k)=2k(\tau-\tau').
\end{equation}
Also, the parity of $\mu_\tau(\gamma)$ is the opposite of the Lefschetz sign
$\epsilon(\gamma)$.

In our two dimensional case, we can explicitly describe the
Conley-Zehnder index and its behavior under multiple covers as
follows.  Fix a nondegenerate periodic orbit $\gamma$ and a trivialization
$\tau\in\mc{T}(\gamma)$.

\begin{proposition}[Conley-Zehnder index of multiple covers]
\label{prop:CZ}
\begin{description}
\item{(a)}
If $\gamma$ is hyperbolic, then there is an integer $n$ such that
\begin{equation}
\label{eqn:CZHyp}
\mu_\tau(\gamma^k)=kn
\end{equation}
for all $k$.  The sign of the eigenvalues $\lambda,\lambda^{-1}$ of
the linearized return map $d\phi^p$ is $(-1)^n$.
\item{(b)}
If $\gamma$ is elliptic, then there is an irrational number $\theta$
such that for all $k$,
\begin{equation}
\label{eqn:CZEll}
\mu_\tau(\gamma^k)=2\lfloor k
\theta\rfloor + 1.
\end{equation}
\end{description}
\end{proposition}

\begin{proof}
This can be seen from the explicit picture of $\op{Sp}(2,\R)$ given by
Segal in \cite{csm}.  For details see \cite{ht1}.  \end{proof}

In case (b), we call $\theta$ the \define{monodromy angle} of
$\gamma$, and we note that the eigenvalues of the linearized return
map are $e^{\pm 2\pi i\theta}$.  Changing the trivialization $\tau$
will change $\theta$ by an integer.  (Perhaps one should refer to
$2\pi\theta$ rather than $\theta$ as the ``angle'', but this would
lead to irritating factors of $2\pi$ later.)

%\begin{remark}
%We can unify the above cases by defining $\theta=n/2$ in the
%hyperbolic case.  Then for any nondegenerate periodic orbit $\gamma$,
%we have
%\[
%\mu_\tau(\gamma^k)=\lfloor k\theta\rfloor + \lceil k\theta\rceil.
%\]
%\end{remark}

\subsection{The relative intersection pairing}
\label{sec:pairing}

\begin{definition}
\label{def:rep}
Let $Z\in H_2(Y;\alpha,\beta)$.  A
\define{representative} of $Z$ is an immersed oriented compact surface
$S$ in
$[0,1]\times Y$ such that:
\begin{itemize}
\item
$\partial S$ consists of positively oriented (resp.\ negatively
oriented) covers of $\{1\}\times\alpha_i$
(resp. $\{0\}\times\beta_j$) whose total multiplicity is $m_i$
(resp. $n_j$).
\item
$[\pi(S)]=Z$, where $\pi:[0,1]\times Y\to Y$ denotes the projection.
\item
$S$ is embedded in $(0,1)\times Y$, and $S$ is transverse to
$\{0,1\}\times Y$.
\end{itemize}
\end{definition}

\begin{definition}
\label{def:tauRep}
If $\tau\in\mc{T}(\alpha,\beta)$, we say that $S$ above is a
\define{$\tau$-representative} of $Z$ if, in addition,
\begin{itemize}
\item
$\pi|_S$ is an immersion near $\partial S$.
\item
$S$ contains $m_i$ (resp. $n_j$) singly covered boundary circles at
$\{1\}\times\alpha_i$ (resp. $\{0\}\times\beta_j$).  The $m_i$
(resp. $n_j$) nonvanishing sections of $E$ over $\alpha_i$
(resp. $\beta_j$), given by projecting the conormal vectors in $S$,
are $\tau$-trivial.  Moreover, in each fiber of $E$ over $\alpha_i$ or
$\beta_j$, these sections lie in distinct rays (emanating from the
origin).
\end{itemize}
\end{definition}

\begin{definition}
\label{def:Q}
If $\tau\in\mc{T}(\alpha,\beta)$ and $Z,Z'\in H_2(Y;\alpha,\beta)$, we
define the \define{relative intersection pairing} $Q_\tau(Z,Z')\in\Z$
as follows.  Let $S$ and $S'$ be $\tau$-representatives of $Z$ and
$Z'$ such that the projected conormal vectors at the boundary all lie
in different rays.  Then $Q_\tau(Z,Z')$ is the signed number of
intersections of $S$ and $S'$ in $(0,1)\times Y$, after perturbing in
the interior to make these intersections transverse.
\end{definition}

We clearly have the symmetry $Q_\tau(Z,Z')=Q_\tau(Z',Z)$.  If we
change one of the relative homology classes or the trivialization, then $Q$
behaves as follows.

\begin{lemma}
\label{lem:pairing}
\begin{description}
\item{(a)} $Q_\tau(Z,Z')$ is well defined, and
\[
Q_\tau(Z_1,Z')-Q_\tau(Z_2,Z')=(Z_1-Z_2)\cap h.
\]
\item{(b)}
$
Q_\tau(Z,Z')-Q_{\tau'}(Z,Z')=\sum_im_i^2({\tau_i'}^+-\tau_i^+) -
\sum_jn_j^2({\tau_j'}^--\tau_j^-).
$
\end{description}
\end{lemma}

Here `$\cap$' indicates homological intersection number in
$Y$.

\begin{proof}
(a) Let $S_1$, $S_2$, and $S'$ be $\tau$-representatives of $Z_1$,
$Z_2$, and $Z'$ respectively.  (To show that $Q_\tau(Z,Z')$ is well
defined, take $Z_1=Z_2=Z$.)  Then $S_1\cup(-S_2)$ is a cycle, which is
homologous to
\[
[pt]\tensor(Z_1-Z_2)\in H_2([0,1]\times Y),
\]
and can be represented by an embedded surface $S''\subset\{1/2\}\times
Y$.  If we perturb $S'$ to be transverse to $\{1/2\}\times Y$ and to
$S''$, then $S'\cap(\{1/2\}\times Y)$ is a 1-manifold in the homology
class $[pt]\tensor h$.  Letting `\#' denote oriented intersection
number, we then have
\[
\begin{split}
Q_\tau(Z_1,Z')-Q_\tau(Z_2,Z') &= \#(S_1\cap S') - \#(S_2\cap S')\\
&= \#(S''\cap S')\\
&= (Z_1-Z_2)\cap h.
\end{split}
\]
In the middle equality above, we have used the fact that the surfaces
are $\tau$-representatives, in order to glue $S_1$ and $S_2$ to a
surface in $(0,1)\times Y$ without introducing or cancelling any
intersections with $S'$ near the boundary.

(b)  If $Z=Z'$, then this follows from
equations \eqref{eqn:writheTriv} and
\eqref{eqn:adjunctionB}, which are proved in
\S\ref{sec:adjunction}. 
One also needs to use the fact that, as in equation
\eqref{eqn:chernTriv},
\[
c_1(N,\tau)-c_1(N,\tau')
=\sum_im_i({\tau_i'}^+-\tau_i^+)-\sum_jn_j({\tau_j'}^--\tau_j^-).
\]
To deduce the general case use (a).  (One can also prove (b) more
directly as in Lemma~\ref{lem:adjunctionC}.)
\end{proof}

\section{Relative adjunction formulas}
\label{sec:adjunction}

We now prove relative adjunction formulas for surfaces with boundary
in $[0,1]\times Y$ or with ends in $\R\times Y$.  These formulas
resemble the adjunction formulas in the closed case, but there is a
correction term involving the writhes of braids determined by the ends
of the surface.  These formulas will play an important role in the
subsequent calculations.

\subsection{Braids and writhe}
\label{sec:braidsIntro}

We begin with some preliminaries on braids.  Let $\gamma$ be a closed
orbit and $\tau\in\mc{T}(\gamma)$.  Let $\xi\subset Y$ be a {\em
braid} in a tubular neighborhood $U$ of $\gamma$.  (That is, $\xi$ is
a compact embedded 1-manifold, and the projection of $\xi$ to $\gamma$
is a submersion.)   Let $A$ denote an annulus, and choose an
identification
\begin{equation}
\label{eqn:annulus}
U\simeq S^1\times\R^2\simeq A\times(0,1)
\end{equation}
whose normal derivative along $\gamma$ agrees with the trivialization
$\tau$. We define the \define{writhe} $w_\tau(\xi)\in\Z$ to be the
signed number of crossings in the projection of $\xi$ to $A$.  Our
sign convention is that counterclockwise twists contribute positively
to the writhe.  (This convention is opposite from much of the knot
theory literature, but works well in this paper, especially in
\S\ref{sec:braids} when we relate writhe to winding numbers.)  
The writhe depends only on the isotopy class of $\xi$ and the homotopy
class of the trivialization $\tau$.  If we change the trivialization,
then
\begin{equation}
\label{eqn:writheTriv}
w_\tau(\xi)-w_{\tau'}(\xi)=m(m-1)(\tau-\tau')
\end{equation}
where $m$ denotes the number of strands in the braid $\xi$.  The reason is
that changing the trivialization by one is equivalent to inserting a
full counterclockwise twist into the braid, which has writhe $m(m-1)$,
as each of the $m$ strands crosses over each of the other $m-1$ strands.

Now let $\alpha$ and $\beta$ be orbit sets and let
$S\subset[0,1]\times Y$ be a representative of $Z\in
H_2(Y;\alpha,\beta)$, as in Definition~\ref{def:rep}.  For $s$ close
to 1, the intersection of $S$ with $\{s\}\times Y$ consists of a braid
$\xi_i^+$ with $m_i$ strands in a neighborhood of $\alpha_i$, for each
$i$. Likewise, for $s$ close to 0 we obtain braids $\xi_j^-$ near
$\beta_j$.  To simplify the notation, if
$\tau\in\mc{T}(\alpha,\beta)$, define the total writhe
\[
w_\tau(S) = \sum_iw_{\tau_i^+}(\xi_i^+) - \sum_jw_{\tau_j^-}(\xi_j^-).
\]

\subsection{Relative adjunction formulas}

If $C$ is a flow line from $\alpha$ to $\beta$ without trivial
cylinders, then under a diffeomorphism of $\R$ with $(0,1)$, we
obtain a surface in $(0,1)\times Y$, whose closure in $[0,1]\times Y$
is a representative of $[C]$ as in Definition~\ref{def:rep}.

\begin{proposition}[relative adjunction formulas]
\label{prop:adjunction}
Let $\alpha,\beta$ be orbit sets, let $S$ be a representative of $Z\in
H_2(Y;\alpha,\beta)$, and let $\tau\in\mc{T}(\alpha,\beta)$.  Let $N$
denote the normal bundle to $S$.  Then:
\begin{description}
\item{(a)}
If $S$ comes from a flow line without trivial cylinders as above, then
\begin{equation}
\label{eqn:adjunctionA}
c_1(E|_Z,\tau) = \chi(S) + c_1(N,\tau).
\end{equation}
\item{(b)}
For general representatives $S$, equation \eqref{eqn:adjunctionA}
holds mod 2, and
\begin{equation}
\label{eqn:adjunctionB}
c_1(N,\tau) = w_\tau(S) + Q_\tau(Z,Z).
\end{equation}
\end{description}
\end{proposition}

\begin{proof}
First note that $c_1(N,\tau)$ is well defined, because projection induces a
canonical isomorphism $N|_{\partial S}= E|_{\partial S}$, so $\tau$
induces a trivialization of $N|_{\partial S}$.

(a) If $S$ comes from a flow line, then up to a homotopy
$(0,1)\approx\R$, we have an isomorphism of complex vector bundles
\begin{equation}
\label{eqn:complexIsomorphism}
(\underline{\C}\oplus E|_S)  \simeq  TS\oplus N,
\end{equation}
since both are isomorphic to $T(\R\times Y)|_S$.  Let $\psi_E$ and
$\psi_N$ be $\tau$-trivial sections of $E|_{\partial S}$ and
$N|_{\partial S}$, and let $\psi_S$ be a nonvanishing section of
$TS|_{\partial S}$ tangent to $\partial S$.  Over $\partial S$, under
the isomorphism
\eqref{eqn:complexIsomorphism}, we have a homotopy through
nonvanishing sections of the determinant line bundles
\begin{equation}
\label{eqn:homotopicSections}
1\wedge\psi_E\approx\psi_S\wedge\psi_N.
\end{equation}
Now in general, if $L_i$ is a complex line bundle on $S$ and $s_i$ is a
nonvanishing section of $L_i|_{\partial S}$ up to homotopy for $i=1,2$, then
\[
c_1(\det(L_1\oplus L_2),s_1\wedge s_2)=c_1(L_1,s_1)+c_1(L_2,s_2).
\]
Applying this identity to both sides of
\eqref{eqn:complexIsomorphism}, with respect to the sections in
\eqref{eqn:homotopicSections}, we obtain \eqref{eqn:adjunctionA}.  

(b)
If $S$ does not come from a pseudoholomorphic curve, then the
isomorphism \eqref{eqn:complexIsomorphism} still holds at the level of
real vector bundles, and still respects the complex structures over
$\partial S$ after straightening $S$ to be normal to
$\{0,1\}\times Y$.  It follows that the relative first Chern
classes differ by an even integer, because changing the complex
structure on a rank 2 complex vector bundle over a closed surface
changes the first Chern class by an even integer.  Thus
\eqref{eqn:adjunctionA} holds mod 2.

To prove \eqref{eqn:adjunctionB}, let $\epsilon>0$ be small and let
$S_0=S\cap ([\epsilon,1-\epsilon]\times Y)$.  Let $S'$ be a surface in
which $S\setminus S_0$ is replaced by a surface $S_1$ consisting of
cobordisms with $\tau$-trivial braids, so that $S'$ is a
$\tau$-representative of $Z$ as in Definition~\ref{def:tauRep}.  Let
$\psi$ be a section of the normal bundle $N'$ to $S'$ that is
$\tau$-trivial over $\partial S_1$.  Let $\psi_0,\psi_1$ denote the
restrictions of $\psi$ to $S_0,S_1$.  We can compute $Q_\tau(Z,Z)$ by
counting the intersections of $S'$ with a pushoff of $S'$ via $\psi$,
so
\begin{equation}
\label{eqn:zeroSection}
\begin{split}
Q_\tau(Z,Z) &= \#\psi^{-1}(0)\\
&= \#\psi_0^{-1}(0) + \#\psi_1^{-1}(0)\\
&= c_1(N,\tau) + \#\psi_1^{-1}(0),
\end{split}
\end{equation}
where `\#' denotes the number of points with signs, after perturbing
to obtain transversality.  Now
\begin{equation}
\label{eqn:correctionTerm}
\#\psi_1^{-1}(0)=-w_\tau(S).
\end{equation}
The reason is that in our cobordism of braids, we can take $\psi_1$ to
be the projection of a nonzero vertical tangent vector
in $A\times (0,1)$.  This section will have
zeroes at the branch points of the projection to
$([0,\epsilon]\cup[1-\epsilon,1])\times Y$, where the writhes of the
braids change.  After an orientation check, this proves
\eqref{eqn:correctionTerm}.  Together with
\eqref{eqn:zeroSection}, this proves \eqref{eqn:adjunctionB}.
\end{proof}

\begin{remark}[adjunction with singularities]
Suppose $C\to\R\times Y$ is a pseudoholomorphic curve which is almost
a flow line, except that it fails to be embedded at a finite number of
singular points.  Let $Z=[C]$.  Then
\begin{equation}
\label{eqn:singular}
c_1(E|_Z,\tau) = \chi(C) +  w_\tau(C) + Q_\tau(Z,Z)- 2\delta(C).
\end{equation}
Here $\delta(C)$ is a sum of positive integer contributions from each
singularity as in \cite{mcduff}.  Namely, near each singular point, we
can perturb the surface to an immersion which is symplectic with
respect to the symplectic form $\omega+ds\wedge dt$ on $\R\times Y$,
and which has only transverse double point singularities; the local
contribution to $\delta$ is then the number of double points.  To
prove \eqref{eqn:singular}, we can carry out the above perturbation
near each singularity without affecting any of the terms in
\eqref{eqn:singular}; then the normal bundle $N\to C$ is well defined,
and a straightforward modification of the proof of
Proposition~\ref{prop:adjunction} shows that
\eqref{eqn:adjunctionA} still holds, while
\eqref{eqn:adjunctionB} holds with a correction term of $2\delta(C)$.
\end{remark}

\subsection{Properties of the relative index}
\label{sec:relativeIndexProperties}

We now prove Proposition~\ref{prop:relativeIndex}.

Part (a) follows from equations \eqref{eqn:chernTriv} and
\eqref{eqn:CZTriv} and
Lemma~\ref{lem:pairing}(b).  The key calculation is that
$\sum_{k=1}^m2k=m^2+m$.

To prove part (b), we compute $\murel(\alpha,\beta;Z)$ and
$\murel(\beta,\gamma;W)$ using trivializations
$\tau_1\in\mc{T}(\alpha,\beta)$ and $\tau_2\in\mc{T}(\beta,\gamma)$
which agree over the orbits in $\beta$; and we compute
$\murel(\alpha,\gamma;Z+W)$ using the corresponding trivialization
$\tau_3\in\mc{T}(\alpha,\gamma)$.  Then
\begin{align*}
c_1(E|_{Z+W},\tau_3)&=
c_1(E|_Z,\tau_1)+c_1(E|_W,\tau_2),\\
Q_{\tau_3}(Z+W,Z+W)&=Q_{\tau_1}(Z,Z)+Q_{\tau_2}(W,W),
\end{align*}
and the Conley-Zehnder index terms over $\beta$ cancel.

To prove part (c), we first observe from Proposition~\ref{prop:CZ} that if
$\alpha$ and $\beta$ are admissible, then
\[
\sum_i\sum_{k=1}^{m_i}\mu_\tau({\alpha_i}^k)
-\sum_j\sum_{k=1}^{n_j}\mu_\tau({\beta_j}^k) \equiv
\sum_i(\epsilon(\alpha_i)+m_i) - \sum_j(\epsilon(\beta_j)+n_j) \mod 2.
\]
Now let $S$ be a $\tau$-representative of $Z$.  Then $w_\tau(S)=0$, so
by the relative adjunction formulas \eqref{eqn:adjunctionA} and
\eqref{eqn:adjunctionB},
\begin{equation}
\label{eqn:parity1}
c_1(E|_Z,\tau)+Q_\tau(Z,Z)\equiv \chi(S) \mod 2.
\end{equation}
By the classification of surfaces, we have
\begin{equation}
\label{eqn:parity2}
\chi(S) \equiv \sum_im_i-\sum_jn_j\mod 2.
\end{equation}
Combining the above three equations completes the proof of part (c).

Part (d) follows from equation \eqref{eqn:chernHom} and
Lemma~\ref{lem:pairing}(a).
\qed

\section{Partitions at the ends of flow lines}
\label{sec:combinatorics}

As remarked in the introduction, a flow line from $\{(\alpha_i,m_i)\}$
to $\{(\beta_j,n_j)\}$ without trivial cylinders determines
partitions of the integers $m_i$ and $n_j$.  In terms of these
partitions, we will now define the notion of ``admissible flow
lines'', for which the dimension of the moduli space can be maximized.

\subsection{Incoming and outgoing partitions}
\label{sec:partitions}

We begin with the definition of some curious partitions.  Let $\theta$
be an irrational number.  If $q$ is a positive integer, define
\[
f_\theta(q)= q^{-1}\lceil q\theta\rceil.
\]
Let $S_\theta$ denote the set of positive integers $q$ such that
$f_\theta(q')>f_\theta(q)$ for all positive integers $q'<q$.
Equivalently, $f_\theta(q)$ is the smallest rational number larger
than $\theta$ with denominator $q$; and $q\in S_\theta$ when $\theta$
is approximated from above by a rational number with denominator $q$,
more closely than any rational number larger than $\theta$ with
denominator less than $q$.

\begin{definition}
If $m$ is a positive integer, the
\define{incoming partition} $\pin(\theta,m)$ is the
partition of $m$ defined inductively as follows: let
$r=\max(S_\theta\cap[1,m])$, and define
\[
\pin(\theta,m)=r\cup\pin(\theta,m-r).
\]
\end{definition}
Note that the use of the symbol `$\cup$' above is a slight abuse of
notation, as a partition may contain the same integer several times.

\begin{remark}
\label{remark:dependence}
The incoming partition depends only on the class $\theta\mod 1$.  If
we vary $\theta$, then $\pin(\theta,m)$ changes only when $\theta$
crosses a rational number with denominator $\le m$.  These facts are
clear from the definition.
\end{remark}

\begin{figure}
\[
\begin{array}{|c||c|c|c|c|c|c|c|}
\hline
 & 2 & 3 & 4 & 5 & 6 & 7 & 8\\ \hline \hline
0,1/8 &  &  &  &  &  &  & 8 \\
\cline{8-8}
1/8,1/7 &  &  &  &  & 6 & \rb{7} & 7,1 \\
\cline{7-8}
1/7,1/6 &  &  & 4 & \rb{5} &  & 6,1 & 6,2 \\
\cline{6-8}
1/6,1/5 &  & 3 &  &  & 5,1 & 5,2 & 5,3 \\
\cline{5-8}
1/5,1/4 &  &  &  & 4,1 & 4,2 & 4,3 & 4,4 \\
\cline{4-8}
1/4,2/7 & 2 &  &  &  &  & 7 & 7,1 \\
\cline{7-8}
2/7,1/3 &  &  & \rb{3,1} & \rb{3,2} & \rb{3,3} & 3,3,1 & 3,3,2 \\
\cline{3-8}
1/3,3/8 &  &  &  &  &  &  & 8 \\
\cline{8-8}
3/8,2/5 &  &  &  & \rb{5} & \rb{5,1} & \rb{5,2} & 5,2,1 \\
\cline{5-8}
2/5,3/7 &  & \rb{2,1} & \rb{2,2} &  &  & 7 & 7,1 \\
\cline{7-8}
3/7,1/2 &  &  &  & \rb{2,2,1} & \rb{2,2,2} & 2,2,2,1 & 2,2,2,2 \\
\cline{2-8}
1/2,4/7 &  &  &  &  & & 7 & 7,1 \\
\cline{7-8}
4/7,3/5 &  &  &  & \rb{5} & \rb{5,1} & 5,1,1 & 5,3 \\
\cline{5-8}
3/5,5/8 &  & \rb{3} & \rb{3,1} &  &  &  & 8 \\
\cline{8-8}
5/8,2/3 &  &  &  & \rb{3,1,1} & \rb{3,3} & \rb{3,3,1} & 3,3,1,1 \\
\cline{3-8}
2/3,5/7 &  &  &  &  &  & 7 & 7,1 \\
\cline{7-8}
5/7,3/4 & 1,1 &  & \rb{4} & \rb{4,1} & \rb{4,1,1} & 4,1,1,1 & 4,4 \\
\cline{4-8}
3/4,4/5 &  &  &  & 5 & 5,1 & 5,1,1 & 5,1,1,1 \\
\cline{5-8}
4/5,5/6 &  & 1,1,1 &  &  & 6 & 6,1 & 6,1,1 \\
\cline{6-8}
5/6,6/7 &  &  & 1,1,1,1 &  &
& 7 & 7,1 \\
\cline{7-8}
6/7,7/8 &  &  &  & \rb{1,\ldots,1} & 1,\ldots,1 &  & 8 \\
\cline{8-8}
7/8,1 &  &  &  &  &  &
\rb{1,\ldots,1} & 1,\ldots,1 \\ \hline
\end{array}
\]
\caption{The incoming partitions
for $2\le m\le 8$ and all $\theta$.  The left column shows the
interval in which $\theta\mod 1$ lies, and the top row indicates $m$.}
\end{figure}

\begin{definition}
We define the \define{outgoing partition}
\[
\pout(\theta,m)=\pin(-\theta,m).
\]
\end{definition}

\begin{remark}
If $m>1$, then $\pin(\theta,m)\cap\pout(\theta,m)=\emptyset$.  (This
makes the gluing theory for periodic Floer homology interesting, see
\cite{ht1}.)
\end{remark}

\begin{proof}
In the following, if $\xi$ is a real number, let $\{\xi\}\in[0,1)$
denote its fractional part.  Fix $q>1$.  If
$\{(q-1)\theta\}>\{q\theta\}$, then a calculation shows that
$f_\theta(q-1)\le f_\theta(q)$, so $q\notin S_\theta$.  This will
happen for $\theta$ or $-\theta$, so
\begin{equation}
\label{eqn:disjoint3}
S_\theta\cap S_{-\theta}=\{1\}.
\end{equation}
One can further show, by induction on the size of the incoming
partition, that if $m>1$, then $\{(m-1)\theta\}>\{m\theta\} \iff
1\in\pin(\theta,m)$.  Applying this observation to $\theta$ and
$-\theta$, we obtain
\begin{equation}
\label{eqn:disjoint4}
1\in\pin(\theta,m)\iff 1\notin\pout(\theta,m).
\end{equation}
Facts \eqref{eqn:disjoint3} and \eqref{eqn:disjoint4} imply the
remark.
\end{proof}

The following properties of the incoming partition will be used in
\S\ref{sec:cylinders} and \S\ref{sec:compactness}.  We leave the
proofs as elementary exercises.

\begin{lemma}
\label{lem:exercise}
Let $\theta\in\R\setminus\Q$ and suppose
$\pin(\theta,m)=(q_1,\ldots,q_k)$.  Then:
\begin{description}
\item{(a)}
$
\sum_{i=1}^k\lceil q_i\theta\rceil = \lceil m\theta\rceil.
$
\item{(b)}
If $i\neq j$, then $\lfloor q_i\theta\rfloor +
\lfloor q_j\theta\rfloor < \lfloor (q_i+q_j)\theta\rfloor$.
\item{(c)}
If $a$ and $b$ are positive integers with
$a+b=q_i$, then $\lfloor a\theta\rfloor + \lfloor b\theta\rfloor =
\lfloor q_i\theta\rfloor$.
\end{description}
\end{lemma}

% I did write proofs of these, but I commented them out:
%
%\begin{proof}
%(a)
%By the definition of $\pin$, we may assume, as in
%\eqref{eqn:decreasing}, that
%\begin{equation}
%\label{eqn:ordConv}
%f_\theta(q_k)\le f_\theta(q')
%\end{equation}
%for all $q'\in\{1,\ldots,m\}$.  Then
%$\pin(\theta,m-q_k)=\{q_1,\ldots,q_{k-1}\}$.  By induction, we may
%assume that
%\[
%\sum_{i=1}^{k-1}\lceil q_i\theta\rceil = \lceil(m-q_k)\theta\rceil.
%\]
%We then need to show that
%$\lceil (m-q_k)\theta\rceil + \lceil q_k\theta\rceil = \lceil
%m\theta\rceil$.
%Equivalently,
%\[
%\{q_k\theta\} \ge \{m \theta\}.
%\]
%Suppose to the contrary that $\{q_k\theta\} < \{m\theta\}$.  This is
%equivalent to
%\[
%\lceil q_k\theta\rceil - q_k\theta > \lceil m\theta\rceil - m\theta.
%\]
%Since $q_k\le m$, we can divide the left side by $q_k$ and the right
%side by $m$, to obtain
%\[
%q_k^{-1}\lceil q_k\theta\rceil > m^{-1}\lceil m\theta\rceil.
%\]
%This contradicts \eqref{eqn:ordConv}, completing the proof.
%
%(b)
%If not, then $\lceil q_1\theta\rceil + \lceil q_2\theta\rceil >
%\lceil (q_1+q_2)\theta\rceil$, and this contradicts
%(a).
%
%(c)
%Since $q_i\in S_\theta$ and $a,b<q_i$, we
%have $q_i^{-1}\lceil q_i\theta\rceil < a^{-1}\lceil
%a\theta\rceil,\,b^{-1}\lceil b\theta\rceil$.  Since $a/q_i+b/q_i=1$, it
%follows that
%\[
%q_i^{-1}\lceil q_i\theta\rceil < \frac{a}{q_i}a^{-1}\lceil a\theta\rceil +
%\frac{b}{q_i} b^{-1}\lceil b\theta\rceil.
%\]
%Multiplying by $q_i$, we obtain $\lceil q_i\theta\rceil < \lceil
%a\theta\rceil + \lceil b\theta\rceil$, which implies (c).
%\end{proof}

\begin{definition}
Let $\gamma$ be a periodic orbit and $m$ a positive integer.  We
define two partitions of $m$, the \define{incoming partition}
$\pin(\gamma,m)$ and the
\define{outgoing partition} $\pout(\gamma,m)$, as follows.
\begin{itemize}
\item
If $\gamma$ is hyperbolic with positive eigenvalues, then
\[
\pin(\gamma,m)=\pout(\gamma,m)=\{1,\ldots,1\}.
\]
\item
If $\gamma$ is hyperbolic with negative eigenvalues, then
\[
\pin(\gamma,m)=\pout(\gamma,m)=\left\{\begin{array}{ll}
\{2,\ldots,2\} & \mbox{if $m$ is even,}\\
\{2,\ldots,2,1\} & \mbox{if $m$ is odd.}
\end{array}\right.
\]
\item If $\gamma$ is elliptic with monodromy angle $\theta$
(see \S\ref{sec:CZ}), then $\pin(\gamma,m)=\pin(\theta,m)$ and
$\pout(\gamma,m)=\pout(\theta,m)$.
\end{itemize}
\end{definition}

\subsection{Admissible flow lines}
\label{sec:admissible}

Let $C$ be a flow line from $\{(\alpha_i,m_i)\}$ to
$\{(\beta_j,n_j)\}$.  Let $C'$ denote the flow line obtained by
deleting all trivial cylinders from $C$, and let $m_i'$ and $n_j'$
denote the total multiplicities of the ends of $C'$ at $\alpha_i$ and
$\beta_j$.  The multiplicities of the outgoing ends of $C'$ at
$\alpha_i$ determine a partition of the integer $m_i'$, which we
denote by $p^+_{i}$.  Likewise the incoming ends of $C'$ determine a
partition $p^-_{j}$ of $n_j'$.  The following definition says that $C$
is admissible if these partitions agree with the incoming and outgoing
partitions defined in the last section, and if the multiplicities of
the trivial cylinders satisfy some additional restrictions.

\begin{definition}
\label{def:admissible}
The flow line $C$ is \define{admissible} if:
\begin{itemize}
\item
For each $i$,
\[
p^+_{i}=p_{\op{out}}(\alpha_i,m_i'),
\]
and for each $k$ with $1\le k\le m_i-m_i'$,
\begin{equation}
\label{eqn:admissible1}
p_{\op{out}}(\alpha_i,m_i'+k)=\pout(\alpha_i,m_i')\cup\pout(\alpha_i,k).
\end{equation}
\item
For each $j$,
\[
p^-_{j}=p_{\op{in}}(\beta_j,n_j'),
\]
and for each $k$ with $1\le k\le n_j-n_j'$,
\begin{equation}
\label{eqn:admissible2}
\pin(\beta_j,n_j'+k)=\pin(\beta_j,n_j')\cup\pin(\beta_j,k).
\end{equation}
\end{itemize}
\end{definition}

\begin{remark}
\label{remark:admissible}
The trivial cylinder conditions \eqref{eqn:admissible1} and
\eqref{eqn:admissible2} hold automatically for hyperbolic orbits with
positive eigenvalues.  If $\alpha_i$ is hyperbolic with negative
eigenvalues, then \eqref{eqn:admissible1} asserts that $m_i=m_i'$ when
$m_i'$ is odd.
\end{remark}

\section{Index theory}
\label{sec:beginComputation}

We now begin to compute the dimension of the moduli space of
flow lines and prove the index theorem \ref{thm:index}.

\subsection{A general index formula}

We first recall a general index formula for certain $\dbar$ operators
on punctured Riemann surfaces.  This formula is proved by Schwarz
\cite{s}, in different notation.

Let $C$ be a Riemann surface of genus $g$ with $k$ punctures.  Let
$V\to C$ be a Hermitian vector bundle of rank $n$ with connection.
Let $S$ be a section of $T^{0,1}C\tensor_{\C}\End_\R(V)$.  Let $p>2$.
We are interested in the real linear differential operator
\begin{equation}
\label{eqn:D}
D=\dbar+S:L^p_1(V)\to L^p(T^{0,1}C\tensor_{\C} V).
\end{equation}

To make sense of this, we need some assumptions.  We choose a
holomorphic identification of each of the $k$ ends of $C$ with
$\R^+\times(\R/\Z)$ with coordinates $s,t$, and we choose a
trivialization $\tau_i$ of $V$ over the $i^{th}$ end of $C$.  In terms
of these choices, on the $i^{th}$ end we can write the operator $D$ as
\begin{equation}
\label{eqn:localForm}
\partial_s+J_0\partial_t+S_i(s,t),
\end{equation}
acting on functions $\R^+\times S^1\to\C^n$. Here $S_i(s,t)$ is a real
square matrix of rank $2n$, and $J_0$ denotes the standard complex
structure on $\C^n$.  We assume that the limit matrix
\begin{equation}
\label{eqn:symmetricLimit}
S_i(t)= \lim_{s\to\infty}S_i(s,t)
\end{equation}
exists and is symmetric.  

We choose a metric on $C$ which on the ends is asymptotically
cylindrical in the above coordinates.  This allows us to define the
spaces $L^p_1$ and $L^p$ in the definition \eqref{eqn:D} of $D$ above,
and the limiting assumption \eqref{eqn:symmetricLimit} implies that
$D$ sends $L^p_1$ to $L^p$.

Next, define a path of symplectic matrices
$\Psi_i(t)$ for $t\in[0,1]$ by
\[
\Psi_i(0)=1,\quad\frac{d\Psi_i(t)}{dt}=J_0S_i(t).
\]
We assume that
\begin{equation}
\label{eqn:nondegenerate}
1\notin\op{Spec}(\Psi_i(1)).
\end{equation}
Then the Maslov index of the path of symplectic matrices
$\Psi_i$ is defined, and we denote it by $\mu_i^\tau\in\Z$.
Let
\[
c_1(V,\tau)\in H^2(C,\partial C)=\Z
\]
denote the relative first Chern class of $V$ with respect to the
trivializations $\tau_i$ at the ends.

\begin{theorem}[general index formula]
\label{thm:schwarz}
\cite{s}
Under the assumptions \eqref{eqn:symmetricLimit} and
\eqref{eqn:nondegenerate}, the operator $D$ is Fredholm, and
\[
\op{ind}(D)=n\chi(C)+2c_1(V,\tau)+\sum_{i=1}^k\mu_i^\tau.
\]
\end{theorem}

\subsection{Beginning the index computation}
\label{sec:begin}

We now use the general index formula to begin the computation of the
dimension of the moduli space of flow lines.  Let $C$ be a flow line
from $\{(\alpha_i,m_i)\}$ to $\{(\beta_j,n_j)\}$, and let $d$ denote
the degree of these orbit sets. Assume for now that $C$ contains no
trivial cylinders.  We write the partitions associated to the ends of
$C$ as $p_i^+=(q_{i,1},q_{i,2},\ldots)$ and
$p_j^-=(q'_{j,1},q'_{j,2},\ldots)$.  To simplify notation, if
$\tau\in\mc{T}(\alpha,\beta)$ is a trivialization, define
\begin{align*}
\mu^0_\tau(C) &= \sum_i
\sum_r\mu_{\tau_i^+}\left(\alpha_i^{q_{i,r}}\right) -
\sum_j \sum_r\mu_{\tau_j^-}\left(\beta_j^{q'_{j,r}}\right),\\
\mu_\tau(C) &= \sum_i\sum_{k=1}^{m_i}\mu_{\tau_i^+}(\alpha_i^k)
-\sum_j\sum_{k=1}^{n_j}\mu_{\tau_j^-}(\beta_j^k).
\end{align*}
Also let $c_\tau(C) = c_1(E|_{[C]},\tau)$ and $Q_\tau(C) = Q_\tau([C],[C])$.

\begin{lemma}
\label{lem:indexBeginning}
Assume $(\phi,J)$ is $d$-admissible and $J$ is generic.  If $C$ is a
flow line without trivial cylinders as above, then $\mc{M}_C$ is a
manifold and
\[
\op{dim}(\mc{M}_C)
\le
c_\tau(C)+Q_\tau(C)+w_\tau(C)+\mu_\tau^0(C).
\]
\end{lemma}

\begin{proof}
A deformation of a flow line might not be embedded, because
singularities might appear at the ends.  Thus we need to consider
``generalized flow lines'' (GFL's), see Definition~\ref{def:GFL}.  Any
GFL near $C$ is ``quasi-embedded'', i.e.\ embedded except possibly for
finitely many singular points, see
\S\ref{sec:GFL}. Let $\widehat{\mc{M}}_C$ denote the component of the
moduli space of quasi-embedded GFL's containing $C$.  A standard
transversality calculation, which we defer to
Lemma~\ref{lem:transversality}(b), shows that $\widehat{\mc{M}}_C$ is
a manifold.  Then
\[
\dim(\widehat{\mc{M}}_C)=\op{ind}(D_C),
\]
where $D_C$ denotes the linearized $\dbar$ operator
\begin{equation}
\label{eqn:dbar}
D_{C}:L^p_1(N)\to L^p(T^{0,1}C\tensor_\C N),
\end{equation}
see \cite[\S3]{ms}.  Here $p>2$, and $N$ denotes the normal bundle to
the embedded curve $C$ in $\R\times Y$.  The spaces $L^p_1$ and $L^p$
are defined using an $\R$-invariant metric on $\R\times Y$.  Note that
as in \cite{t2}, we use the normal bundle, rather than the restriction
of the tangent bundle of $\R\times Y$ to $C$, because we are not
fixing the complex structure on the domains of our pseudoholomorphic
curves.

We now compute that
\[
\begin{split}
\op{ind}(D_{C}) &=\chi(C)+2c_1(N,\tau)
+\mu^0_\tau(C)\\
&= c_\tau(C)+c_1(N,\tau)+\mu^0_\tau(C)\\
&= c_\tau(C)+Q_\tau(C)+w_\tau(C)+\mu^0_\tau(C).
\end{split}
\]
Here the first equality follows from Theorem~\ref{thm:schwarz}; the
operator $D_C$ has the local form \eqref{eqn:localForm} at the ends of
$C$, and the Conley-Zehnder indices come out right, by equations
\eqref{eqn:A} and \eqref{eqn:pseudoholomorphic} below.  The second and
third equalities follow from the relative adjunction formulas
\eqref{eqn:adjunctionA} and \eqref{eqn:adjunctionB}.

By equation \eqref{eqn:singular}, any GFL $C'\in\mc{M}_C$ satisfies
 $w_\tau(C')=w_\tau(C)$.  It follows from this and the discussion in
\S\ref{sec:braids} that $\mc{M}_C$ is a component of a stratum of
$\widehat{\mc{M}}_C$ in which certain coefficients in asymptotic
expansions of the ends vanish.  Applying the same transversality
calculation as in Lemma~\ref{lem:transversality}(b) to GFL's
satisfying this restriction, we find that $\mc{M}_C$ is a submanifold
of $\widehat{\mc{M}}_C$ for generic $J$.  (But we remark that in the
important case when $w_\tau(C)$ saturates the upper bound
\eqref{eqn:braids} below, equation
\eqref{eqn:singular} implies that in fact
$\mc{M}_C=\widehat{\mc{M}}_C$.)
\end{proof}

\section{Braids at the ends of flow lines}
\label{sec:braids}

This section is devoted to proving the following proposition.  We use
the notation of \S\ref{sec:begin}.

\begin{proposition}[braids]
\label{prop:braids}
If $C$ is a flow line without trivial cylinders, and if $(\phi,J)$ is
admissible near the ends of $C$, then
\begin{equation}
\label{eqn:braids}
w_\tau(C)+\mu_\tau^0(C)\le\mu_\tau(C).
\end{equation}
Equality holds only if $C$ is admissible.
\end{proposition}

The idea of the proof is to show that the braids at the ends of $C$ are
iterated nested cablings of torus braids, and to bound their writhe
using bounds on winding numbers in terms of the Maslov index from
\cite{hwz2}.  Proposition~\ref{prop:braids}, together with
Lemma~\ref{lem:indexBeginning}, implies Theorem~\ref{thm:index} for
flow lines without trivial cylinders.  We will deal with trivial
cylinders in \S\ref{sec:cylinders}.

\subsection{Asymptotics and winding numbers}

We begin by describing the basic structure of an incoming end of
multiplicity $q$ at a periodic orbit $\gamma$.  Similar results hold
for outgoing ends, but by symmetry we do not need to consider
these, see Lemma~\ref{lem:outgoingInequality}.

\begin{assumption}
For the rest of \S\ref{sec:braids}, we assume that $(\phi,J)$ is
admissible near $\gamma$.
\end{assumption}
So we can choose a neighborhood $N$ of $\gamma$
and an identification
\begin{equation}
\label{eqn:linear}
N\simeq S^1\times\R^2,
\end{equation}
identifying $\gamma\simeq S^1\times\{0\}$, compatible with the
projection $N\to\gamma$ induced by the projection $Y\to S^1$, such
that with respect to this identification:
\begin{itemize}
\item
In $N$, parallel transport via the natural connection on $Y\to S^1$ is
linear.
\item
The restriction to $E$ of the almost complex structure $J$ is given by
a constant matrix on each fiber of $N\to S^1$.
\end{itemize}
The normal derivative along $\gamma$ of the identication
\eqref{eqn:linear} determines a trivialization $\tau$ of $E|_\gamma$.  With
respect to the identification \eqref{eqn:linear} (resp. the
trivialization $\tau$), the covariant derivative in $S^1\times\R^2$
near $\gamma$ (resp. $E$) is given by
\[
\nabla_t=\partial_t-J(t)S(t).
\]
Here $t$ denotes the $S^1$ coordinate, $J(t)$ is the restriction to
$E$ of the almost complex structure, and $S(t)$ is a $2\times 2$
matrix which is symmetric with respect to the metric determined by
$J(t)$ and $\omega$.

Suppose a flow line has an incoming end at $\gamma$ of multiplicity
$q$.  For $s<<0$, the projection of the end to $\R\times S^1$ is a
covering without branch points.  Thus for $R<<0$, we can describe the
end by a pseudoholomorphic map of bundles over $(-\infty,R)\times
S^1$,
\[
u:(-\infty,R)\times\tilde{S^1} \longrightarrow (-\infty,R)\times Y,
\]
where $\tilde{S^1}$ denotes a $q$-fold connected covering of $S^1$.

Let $\tilde{t}$ denote the $\tilde{S^1}$ coordinate, and let $t$
denote the projection of $\tilde{t}$ to $S^1$.  Define a self-adjoint
operator
\begin{equation}
\label{eqn:A}
A=-J(t)\nabla_{\tilde{t}}
=
-J(t)\partial_{\tilde{t}}-S(t):
C^\infty(\tilde{S^1},\R^2)\to
C^\infty(\tilde{S^1},\R^2).
\end{equation}
By analytic perturbation theory \cite{k}, as used in \cite{hwz2},
there is a countable set of eigenfunctions $\{e_n\}$ with
$Ae_n=\lambda_ne_n$, which constitute an orthonormal basis for
$L^2(\tilde{S^1},\R^2)$ over $\R$.

\begin{lemma}[asymptotic expansion]
\label{lem:expansion}
If a flow line has an incoming end at $\gamma$ of multiplicity $q$,
then we can expand the end for $s<<0$ as
\begin{equation}
\label{eqn:expansion}
u(s,\tilde{t})=\sum_na_ne^{\lambda_ns}e_n(\tilde{t})
\end{equation}
with $a_n\in\R$.  If $a_n\neq 0$, then $\lambda_n$ is positive.
\end{lemma}

\begin{proof}
For fixed $s$, we can expand $u(s,\cdot)$ as a linear combination of the
$e_n$'s.  By the admissibility conditions at the beginning of this section,
$u$ is pseudoholomorphic near $\gamma$ if and only if it satisfies the
linear equation
\begin{equation}
\label{eqn:pseudoholomorphic}
(\partial_s-A)u(s,\cdot)=0.
\end{equation}
This implies \eqref{eqn:expansion}.  Now we must have
$u(s,\tilde{t})\to 0$ as $s\to-\infty$.  In particular
\[
a_ne^{\lambda_ns}=\langle u(s,\cdot),e_n\rangle\to 0,
\]
which implies that $a_n=0$ unless $\lambda_n>0$.
\end{proof}

If $e$ is a (nonzero) eigenfunction of $A$, then $e(\tilde{t})\neq 0$
for all $\tilde{t}$, by uniqueness of solutions to ODE's.  Thus we can
define the \define{winding number} $\eta(e)\in\Z$ of the path $e$ in
$\R^2$ around the origin.

\begin{lemma}[eigenvalues and winding numbers]
\label{lem:hwz}
\begin{description}
\item{(a)}
If $e,e'$ are eigenfunctions of $A$ corresponding to eigenvalues
$\lambda\le\lambda'$, then $\eta(e)\le\eta(e')$.
\item{(b)}
For each integer $\eta$, the space of eigenfunctions with winding
number $\eta$ is 2-dimensional.
\item{(c)}
If $\gamma$ is nondegenerate (so that $0$ is not an eigenvalue of
$A$), then the maximal winding number for a negative eigenvalue is
$\lfloor\mu_\tau(\gamma^q)/2\rfloor$, and the minimal winding number
for a positive eigenvalue is $\lceil\mu_\tau(\gamma^q)/2\rceil$.
\end{description}
\end{lemma}

\begin{proof}
This is all proved in \cite[\S3]{hwz2}.
\end{proof}

\begin{example}
For the simplest kind of elliptic end with monodromy angle $\theta$,
if $q=1$, $J(t)=\begin{pmatrix}0&-1\\1&0\end{pmatrix}$, and
$J(t)S(t)=\begin{pmatrix}0&-2\pi\theta\\2\pi\theta&0\end{pmatrix}$,
then identifying $\R^2$ with $\C$, we can take $e_n(t)=e^{2\pi int}$
for $n\in\Z$, with eigenvalues $\lambda_n=2\pi(n-\theta)$.
\end{example}

Now let $\xi$ be the braid corresponding to an incoming end at
$\gamma$ of multiplicity $q$, for $s<<0$.  We assume
that this end is \define{nontrivial}, i.e.\ not a trivial cylinder.
We temporarily ignore any other ends of the flow line at $\gamma$.
Let $\eta_\tau(\xi)$ denote the winding number of $\xi$ around
$\gamma$, with respect to the trivialization $\tau$.

\begin{lemma}[winding bound]
\label{lem:windingBound}
The winding number $\eta_\tau(\xi)$ is well defined, and
\begin{equation}
\label{eqn:windingBound}
\eta_\tau(\xi)\ge\left\lceil\frac{\mu_\tau(\gamma^q)}{2}\right\rceil.
\end{equation}
\end{lemma}

\begin{proof}
When $s<<0$, the smallest eigenvalue $\lambda$ in the expansion
\eqref{eqn:expansion} dominates, so
\[
u(s,\tilde{t})\sim e^{\lambda s}\varphi(\tilde{t})
\]
where $A\varphi=\lambda\varphi$.  Thus $\eta_\tau(\xi)=\eta(\varphi)$.
By Lemma~\ref{lem:expansion}, $\lambda>0$, so by
Lemma~\ref{lem:hwz}(c),
$\eta(\varphi)\ge\lceil\mu_\tau(\gamma^q)/2\rceil$.
\end{proof}

\subsection{The writhe of a single end}
\label{sec:writhe}

\begin{lemma}
\label{lem:cabling}
Suppose a flow line has a nontrivial incoming end at $\gamma$ with
multiplicity $q$.  Let $\xi$ be the corresponding braid, ignoring any
other ends at $\gamma$. Then the writhe is bounded by
\begin{equation}
\label{eqn:cabling}
w_\tau(\xi)\ge(q-1)\eta_\tau(\xi).
\end{equation}
\end{lemma}

\begin{proof}
We will prove inductively that \eqref{eqn:cabling} holds for any braid
$\xi$ coming from an expansion \eqref{eqn:expansion} for $s<<0$.
Write $\eta=\eta_\tau(\xi)$.  We begin with the case in which
\begin{equation}
\label{eqn:relativelyPrime}
\op{gcd}(q,\eta)=1.
\end{equation}
Let $\varphi$ denote the dominating eigenfunction in the expansion
\eqref{eqn:expansion}, as in the proof of
Lemma~\ref{lem:windingBound}.  Let $\overline{\xi}$ denote the braid
swept out by $\varphi$.

Claim: over each point in $S^1\times\{0\}\simeq\gamma$, each ray in
$\R^2$ intersects at most one strand of the braid $\overline{\xi}$;
and hence the same is true for $\xi$.  To prove the claim, suppose to
the contrary that for some $\tilde{t}\in\tilde{S^1}\simeq\R/q\Z$ and
some $j\in\{1,\ldots,q-1\}$, there is a positive real number $\rho$
with
\[
\varphi(\tilde{t})=\rho\varphi(\tilde{t}+j).
\]
By uniqueness of solutions to the ODE
\[
[-J(t)\partial_{\tilde{t}}-S(t)]\varphi(\tilde{t})=\lambda\varphi(\tilde{t}),
\]
we would have
\[
\varphi(\tilde{t}+z)=\rho\varphi(\tilde{t}+j+z)
\]
for all $z\in\tilde{S^1}$.  Taking $z=j,2j,\ldots$, it follows that
$\rho=1$ and
\[
\varphi(z)=\varphi(\op{gcd}(j,q)+z)
\]
for all $z\in\tilde{S^1}$.  This periodicity implies that the winding
number $\eta$ is divisible by $q/\op{gcd}(j,q)$.  This contradicts our
assumption \eqref{eqn:relativelyPrime}, since this divisor also
divides $q$ and is larger than 1.  This proves the claim.

From the claim, it follows that $\xi$ is isotopic to a $(q,\eta)$
torus braid, and this braid has writhe $\eta(q-1)$.

Suppose now that
\[
\op{gcd}(q,\eta)=d>1.
\]
One can obtain an eigenfunction on $\tilde{S^1}$ with winding number
$\eta$ by pulling back an eigenfunction on a $(q/d)$-fold covering
$\hat{S^1}$ of $S^1$ with winding number $\eta/d$.  By
Lemma~\ref{lem:hwz}(b), this construction yields all eigenfunctions
with winding number $\eta$.  In particular, the dominating
eigenfunction $\varphi_1$ in $\xi$ is pulled back from an
eigenfunction on $\hat{S^1}$, whose braid $\xi_1$ has $q/d$ strands
and winding number $\eta/d$.  It follows that $\xi$ is the
\define{cabling} of $\xi_1$ by some braid $\xi_2$ with $d$ strands.
That is, $\xi$ is obtained by replacing the string of $\xi_1$ with the
braid $\xi_2$ in a tubular neighborhood of $\xi_1$.  The braid $\xi_2$
is obtained by subtracting the $\varphi_1$ term from the expansion of
$\xi$, and lifting the resulting braid to a $(q/d)$-fold covering of
the tubular neighborhood $N$.  Thus $\xi_2$ is approximated by an
eigenfunction $\varphi_2$ with eigenvalue larger than that of
$\varphi_1$.  Let $\eta'$ denote the winding number of $\xi_2$.  Then
by Lemma~\ref{lem:hwz}(a),
\begin{equation}
\label{eqn:eta'}
\eta'\ge\eta.
\end{equation}

Since $\xi$ is the cabling of $\xi_1$ by $\xi_2$, it follows from the
definition of writhe that
\begin{equation}
\label{eqn:writheOfCabling}
w_\tau(\xi)=d^2w_\tau(\xi_1)+w_\tau(\xi_2).
\end{equation}
By induction on $q$, we may assume that
\begin{equation}
\label{eqn:writheOfPieces}
\begin{split}
w_\tau(\xi_1) &\ge\frac{\eta}{d}\left(\frac{q}{d}-1\right),\\
w_\tau(\xi_2) &\ge\eta'(d-1).
\end{split}
\end{equation}
Putting \eqref{eqn:writheOfPieces} into \eqref{eqn:writheOfCabling}
and using \eqref{eqn:eta'} gives \eqref{eqn:cabling}.
\end{proof}

\begin{lemma}[writhe bound]
\label{lem:writheBound}
Suppose a flow line has a nontrivial incoming end at $\gamma$ with
multiplicity $q$ and braid $\xi$, ignoring all other ends.  Then
\begin{equation}
\label{eqn:writheBound}
w_\tau(\xi)\ge\left\lceil\frac{\mu_\tau(\gamma^q)}{2}\right\rceil(q-1).
\end{equation}
If equality holds, then:
\begin{description}
\item{(i)} 
If $\gamma$ is hyperbolic with positive eigenvalues then $q=1$.
\item{(ii)}
If $\gamma$ is hyperbolic with negative eigenvalues, then $q$ is odd or $q=2$.
\end{description}
\end{lemma}

\begin{proof}
Putting Lemma~\ref{lem:windingBound} into Lemma~\ref{lem:cabling}, we
obtain the inequality \eqref{eqn:writheBound}.  Now suppose that
equality holds in \eqref{eqn:writheBound}.  Since $q>1$, equality must
also hold in \eqref{eqn:windingBound}:
\begin{equation}
\label{eqn:windingEquality}
\eta = \eta_\tau(\xi) = \left\lceil\frac{\mu_\tau(\gamma^q)}{2}\right\rceil.
\end{equation}
Suppose (i) or (ii) fails; we will obtain a contradiction.

\underline{Case 1}:
$\gamma$ is hyperbolic with positive eigenvalues and $q>1$.  By
Proposition~\ref{prop:CZ} and equation \eqref{eqn:CZTriv}, we can choose the
trivialization $\tau$ so that $\mu_\tau(\gamma^q)=0$.  By
\eqref{eqn:windingEquality}, $\eta=0$.  As in the proof of
Lemma~\ref{lem:cabling}, we deduce that $\xi$ is the cabling of a
one-strand braid $\xi_1$ with winding number zero by a $q$-strand
braid $\xi_2$ with winding number $\eta'$, whose dominant eigenvalue
is greater than that of $\xi$.  The latter eigenvalue is positive, and by
Lemma~\ref{lem:hwz}(c) there is an eigenfunction with a negative
eigenvalue and the same winding number.  So by Lemma~\ref{lem:hwz}(b),
$\eta'\neq\eta$, and then by Lemma~\ref{lem:hwz}(a), $\eta'>\eta$.  By
\eqref{eqn:writheOfCabling} and \eqref{eqn:writheOfPieces}, we deduce
that $w_\tau(\xi)\ge \eta'(q-1)>0$.  So equality does not hold in
\eqref{eqn:writheBound}, and this is a contradiction.

\underline{Case 2}:
$\gamma$ is hyperbolic with negative eigenvalues, $q$ is even, and
$q>2$.  By Proposition~\ref{prop:CZ} and equation \eqref{eqn:CZTriv}, we can
choose the trivialization $\tau$ so that $\mu_\tau(\gamma^q)=q$.  By
\eqref{eqn:windingEquality}, $\eta=q/2$.  Similarly to Case 1, we
deduce that $\xi$ is the cabling of a $2$-strand braid $\xi_1$ by a
$q/2$-strand braid $\xi_2$ with winding number $\eta'>\eta$.  By
\eqref{eqn:writheOfCabling} and \eqref{eqn:writheOfPieces}, and using
the fact that $d=q/2>1$, we obtain $w_\tau(\xi)>(q/2)(q-1)$.  So
equality does not hold in \eqref{eqn:writheBound}, a contradiction.
\end{proof}

\subsection{Linking of two ends}
\label{sec:linking}

Let $\xi_1$ and $\xi_2$ be two disjoint braids in a neighborhood of
$\gamma$, and let $\tau\in\mc{T}(\gamma)$.  We define the
\define{linking number}
\[
\ell_\tau(\xi_1,\xi_2)\in\Z
\]
to be one half the signed number of crossings of a strand of $\xi_1$
with a strand of $\xi_2$, in the projection to $A$ from equation
\eqref{eqn:annulus}, using the same sign convention as for the writhe
in \S\ref{sec:braidsIntro}.  The linking number is clearly symmetric:
$\ell_\tau(\xi_1,\xi_2)=\ell_\tau(\xi_2,\xi_1)$.

\begin{lemma}[linking bound]
\label{lem:linkingBound}
Suppose a flow line has a nontrivial incoming end at $\gamma$ of
multiplicity $q_i$ with braid $\xi_i$ and winding number $\eta_i$ for
$i=1,2$.  Then
\begin{equation}
\label{eqn:linking}
\ell_\tau(\xi_1,\xi_2)\ge\min(q_1\eta_2,q_2\eta_1).
\end{equation}
\end{lemma}

\begin{proof}
Let $A_q$ denote the operator
\eqref{eqn:A} on a $q$-fold cover of $\gamma$.  Let $\lambda_i$ denote
the smallest eigenvalue of $A_{q_i}$ in the expansion of $\xi_i$.  WLOG
$\lambda_1\le\lambda_2$.  We can pull back the corresponding
eigenfunctions to a $q_1q_2$-fold cover of $S^1$, and applying
Lemma~\ref{lem:hwz}(a) to $A_{q_1q_2}$, we find that
\begin{equation}
\label{eqn:monotonic}
\eta_1q_2\le\eta_2q_1.
\end{equation}
If $\lambda_1<\lambda_2$, then for $s<<0$,
equation~\eqref{eqn:expansion} implies that the braid $\xi_2$ is
{\em nested\/} inside $\xi_1$, i.e.\ there is a tube containing $\gamma$
and $\xi_2$ but not intersecting $\xi_1$.  It follows from the
definition of linking number that
\begin{equation}
\label{eqn:eta1q2}
\ell_\tau(\xi_1,\xi_2)=\eta_1q_2.
\end{equation}
Together with \eqref{eqn:monotonic}, this proves \eqref{eqn:linking}.

If $\lambda_1=\lambda_2$, then equation \eqref{eqn:monotonic} gives
$\eta_1q_2=\eta_2q_1$.  If the coefficients $a$ of the corresponding
eigenfunctions in the expansions \eqref{eqn:expansion} of $\xi_1$ and
$\xi_2$ are different, then for $s<<0$, one can isotope one of the
braids radially towards $\gamma$ without intersecting the other, so we
obtain
\eqref{eqn:eta1q2}, which proves \eqref{eqn:linking} as before.

Otherwise, let $\lambda$ be the smallest eigenvalue of $A_{q_1q_2}$
for which the coefficients of the corresponding eigenfunctions in
$\xi_1$ and $\xi_2$ are different.  Let $\xi_3$ be the braid obtained
from the parts of the expansions of $\xi_1$ and $\xi_2$ involving
eigenvalues smaller than $\lambda$.  Let $q_3$ be the number of
strands of $\xi_3$, and let $w_3$ denote its writhe in the
trivialization $\tau$.  Then for $s<<0$ and $i=1,2$, the braid $\xi_i$
is a cabling of $\xi_3$ by some braid $\xi_i'$ with $q_i/q_3$ strands,
obtained from the remaining terms in the expansion, and we can perform
both cablings in the same tubular neighborhood of $\xi_3$.  It
follows from the definition of linking number that
\begin{equation}
\label{eqn:linkingOfCablings}
\ell_\tau(\xi_1,\xi_2) =
w_3\frac{q_1}{q_3}\frac{q_2}{q_3} + \ell_\tau(\xi_1',\xi_2').
\end{equation}
As in the proof of Lemma~\ref{lem:cabling}, the braid $\xi_i'$ has
winding number $\eta_i'\ge\eta_i$.  By the above reasoning, WLOG,
\begin{equation}
\label{eqn:pullApart}
\ell_\tau(\xi_1',\xi_2')=\eta_1'\frac{q_2}{q_3}\ge\eta_1\frac{q_2}{q_3}.
\end{equation}
Now $\xi_3$ has winding number
\[
\eta_3=\eta_1\frac{q_3}{q_1}=\eta_2\frac{q_3}{q_2},
\]
and the proof of Lemma~\ref{lem:cabling} shows that
\begin{equation}
\label{eqn:w3}
w_3\ge\eta_3(q_3-1)=\frac{\eta_1q_3}{q_1}(q_3-1).
\end{equation}
Putting the inequalities \eqref{eqn:w3} and \eqref{eqn:pullApart} into
equation \eqref{eqn:linkingOfCablings}, we obtain
$\ell_\tau(\xi_1,\xi_2)\ge\eta_1q_2$, which proves \eqref{eqn:linking}.
\end{proof}

\subsection{A combinatorial lemma}

As we will see in \S\ref{sec:noCylinders}, the inequalities we have
just established reduce the proof of Proposition~\ref{prop:braids} to
the following combinatorial lemma.

\begin{lemma}[workhorse inequality]
\label{lem:workhorse}
Let $\gamma$ be a periodic orbit, let $\tau\in\mc{T}(\gamma)$, let
$q_1,\ldots,q_k$ be positive integers, $n=\sum_iq_i$, and
$\rho_i=\left\lceil \mu_\tau(\gamma^{q_i}) / 2 \right\rceil$. Then
\begin{equation}
\label{eqn:workhorse}
\sum_{i=1}^k(\mu_\tau(\gamma^{q_i})-\rho_i)+
\sum_{i,j=1}^k\min(q_i\rho_j,q_j\rho_i)
\ge
\sum_{i=1}^n\mu_\tau(\gamma^i).
\end{equation}
Equality holds if and only if:
\begin{description}
\item{(i)}
If $\gamma$ is hyperbolic with negative eigenvalues, then all $q_i$'s
are even, except that one $q_i$ might equal one.
\item{(ii)}
If $\gamma$ is elliptic with monodromy angle $\theta$, then
\[
\{q_1,\ldots,q_k\}=\pin(\theta,n).
\]
\end{description}
\end{lemma}

\begin{proof}  We first note that the validity of the inequality
\eqref{eqn:workhorse} does not depend on the choice of trivialization,
as changing $\tau$ changes both sides of \eqref{eqn:workhorse}
equally, by adding an integer multiple of $n^2+n$.  We now consider
three cases.

\underline{Case 1}: $\gamma$ is hyperbolic with positive eigenvalues.  We can
choose the trivialization $\tau$ so that $\mu_\tau(\gamma^i)=0$ for
all $i$.  Then the inequality \eqref{eqn:workhorse} trivially holds,
and is an equality, as all terms in it are equal to zero.

\underline{Case 2}: $\gamma$ is elliptic with monodromy angle $\theta$ with
respect to the trivialization $\tau$.  Then $\rho_i=\lceil
q_i\theta\rceil$.  Define
\[
M_\theta(q_1,\ldots,q_k) = 
\sum_{i=1}^k\lfloor q_i\theta\rfloor +
\sum_{i,j=1}^k
\min\left(q_i\lceil q_j\theta\rceil,q_j\lceil q_i\theta\rceil\right)
-n-2\sum_{i=1}^n\lfloor i\theta\rfloor.
\]
Lemma~\ref{lem:workhorse} reduces in this case to the following lemma.
This characterization of the incoming partition will also be useful in
\S\ref{sec:cylinders}.

\begin{lemma}[incoming partitions]
\label{lem:incoming}
If $\theta$ is an irrational number, and $q_1,\ldots,q_k$ are positive
integers whose sum is $n$, then:
\begin{description}
\item{(a)} $M_\theta(q_1,\ldots,q_k)\ge 0$.
\item{(b)} Equality holds in (a) if and only if
$\{q_1,\ldots,q_k\}=\pin(\theta,n)$.
\end{description}
\end{lemma}

\begin{proof}
Recall the notation $f_\theta(q)=q^{-1}\lceil q\theta\rceil$ from
\S\ref{sec:partitions}.  By symmetry, we may assume that
\begin{equation}
\label{eqn:decreasing}
f_\theta(q_1)\ge f_\theta(q_2) \ge \cdots \ge f_\theta(q_k).
\end{equation}
(a) By induction, it will suffice to show that
$M_\theta(q_1,\ldots,q_k)\ge M_\theta(q_1,\ldots,q_{k-1})$.  Under the
ordering convention \eqref{eqn:decreasing}, the latter inequality becomes
\[
\lfloor q_k\theta\rfloor + q_k\lceil q_k\theta\rceil +
2n'\lceil q_k\theta\rceil \ge q_k + 2\sum_{i=n'+1}^n\lfloor i\theta\rfloor
\]
where $n'=q_1+\cdots+q_{k-1}$.  Since $\theta$ is irrational, we can
convert floors to ceilings to rewrite this as $z(\theta)\ge 0$, where
\begin{equation}
\label{eqn:z}
z(\xi)=(2n'+q_k+1)\lceil q_k\xi\rceil + q_k-1 - 2\sum_{i=n'+1}^n\lceil
i\xi\rceil.
\end{equation}
Let $\theta'=f_\theta(q_k)$.  We observe that
\begin{equation}
\label{eqn:zMonotone}
z(\theta)\ge z(\theta'),
\end{equation}
because as a function of $\xi\in[\theta,\theta']$, the first term on
the right side of \eqref{eqn:z} is constant, while the remaining terms
are monotone decreasing.  So to prove (a), it will suffice to show that
\begin{equation}
\label{eqn:ztheta'}
z(\theta')\ge 0.
\end{equation}
Let $\epsilon_i(\xi)=\lceil i\xi\rceil - i\xi$.  Then \eqref{eqn:z}
can be rewritten as
\begin{equation}
\label{eqn:z2}
z(\xi) = (2n'+q_k+1)\epsilon_{q_k}(\xi)+q_k-1-2\sum_{i=n'+1}^n\epsilon_i(\xi).
\end{equation}
The numbers $\epsilon_{n'+1}(\theta'),\ldots,\epsilon_n(\theta')$ are
evenly spaced around the circle $\R/\Z$ with spacing $1/s$, where
$s=q_k/\op{gcd}(\lceil q_k\theta\rceil,q_k)$. So
\begin{equation}
\label{eqn:s}
\begin{split}
\sum_{i=n'+1}^n\epsilon_i(\theta')
&= \frac{q_k}{s}\left(\frac{1}{s}+\frac{2}{s}+\cdots+\frac{s-1}{s}\right)\\
&=\frac{q_k}{2}\left(\frac{s-1}{s}\right)\\ &\le\frac{q_k-1}{2}.
\end{split}
\end{equation}
Putting this inequality into equation \eqref{eqn:z2} with
$\xi=\theta'$, and noting that $\epsilon_{q_k}(\theta')=0$, we obtain
the desired inequality
\eqref{eqn:ztheta'}, so (a) is proved.

(b) We have $M_\theta(q_1,\ldots,q_k)=0$ if and only if equality holds
in the inequalities \eqref{eqn:zMonotone} and \eqref{eqn:s} at each
stage in the above induction.  Now equality holds in
\eqref{eqn:zMonotone} at each stage if and only if
\begin{description}
\item{(i)}
Under the ordering convention \eqref{eqn:decreasing}, the interval
$(\theta,f_\theta(q_i))$ contains no rational numbers of the form
$p/r$ with
\begin{equation}
\label{eqn:range}
q_1+\cdots+q_{i-1}<r\le q_1+\cdots+q_i.
\end{equation}
\end{description}
And equality holds in \eqref{eqn:s} at each stage if and only if
\begin{description}
\item{(ii)}
$\op{gcd}(q_i,\lceil q_i\theta\rceil)=1$.
\end{description}
For a given integer $r$, there exists an integer $p$ with
$p/r\in(\theta,f_\theta(q_i))$ if and only if
$f_\theta(r)<f_\theta(q_i)$.  Thus condition (i) is equivalent to
\begin{description}
\item{(i$'$)}
$f_\theta(q_i)\le f_\theta(r)$ for all $r$ in the range \eqref{eqn:range}.
\end{description}
Also, condition (ii) is equivalent to
\begin{description}
\item{(ii$'$)}
$f_\theta(q_i)\neq f_\theta(r)$ for all $r<q_i$.
\end{description}
This is because if $f_\theta(q_i)=f_\theta(r)$ with $r<q_i$,
i.e. $\lceil q_i\theta\rceil/q_i=\lceil r\theta\rceil/r$, then $q_i$
and $\lceil q_i\theta\rceil$ have a common factor; while conversely,
if $q_i$ and $\lceil q_i\theta\rceil$ are both divisible by $d$, then
$f_\theta(q_i)=f_\theta(q_i/d)$.

Now suppose that conditions (i$'$) and (ii$'$) hold.  By (i$'$) and
the ordering convention \eqref{eqn:decreasing}, $f_\theta(q_i)\le
f_\theta(r)$ for all $r\in[1,q_1+\cdots+q_i]$.  Together with (ii$'$),
this implies that
\begin{equation}
\label{eqn:incoming'}
q_i=\max(S_\theta\cap[1,q_1+\cdots+q_i]).
\end{equation}
Taking $i=k,k-1,\ldots$, we deduce inductively that
$\{q_1,\ldots,q_k\}=\pin(\theta,n)$.  Conversely, if
$\{q_1,\ldots,q_k\}=\pin(\theta,n)$, then by the ordering
\eqref{eqn:decreasing}, equation \eqref{eqn:incoming'} holds by
the definition of $\pin(\theta,n)$, and this implies (i$'$) and (ii$'$).

This completes the proof of Lemma~\ref{lem:incoming}.
\end{proof}

\underline{Case 3} of the proof of Lemma~\ref{lem:workhorse}: $\gamma$ is
hyperbolic with negative eigenvalues.  We can choose the
trivialization $\tau$ so that $\mu_\tau(\gamma^q)=q$ and
$\rho_i=\lceil q_i/2\rceil$.  We then need to show that
$M_{1/2}(q_1,\ldots,q_k)\ge 0$, where
\[
M_{1/2}(q_1,\ldots,q_k) =
-\sum_{i=1}^k\left\lceil\frac{q_i}{2}\right\rceil +
\sum_{i,j=1}^k\min\left(q_i\left\lceil\frac{q_j}{2}\right\rceil,
q_j\left\lceil\frac{q_i}{2}\right\rceil\right) - \frac{n(n-1)}{2},
\]
with equality if and only if all $q_i$'s are even, except that one
$q_i$ might equal 1.  The proof parallels the proof of
Lemma~\ref{lem:incoming}, but is simpler.  Without loss of generality,
\[
q_1^{-1}\left\lceil\frac{q_1}{2}\right\rceil \ge
q_2^{-1}\left\lceil\frac{q_2}{2}\right\rceil \ge \cdots \ge
q_k^{-1}\left\lceil\frac{q_k}{2}\right\rceil.
\]
We then calculate that
\[
M_{1/2}(q_1,\ldots,q_k)-M_{1/2}(q_1,\ldots,q_{k-1}) =
\left(\left\lceil\frac{q_k}{2}\right\rceil-\frac{q_k}{2}\right)(2n-q_k-1).
\]
This is nonnegative since $1\le q_k\le n$, and zero if and
only if $q_k$ is even or $2n-q_k-1=0$, i.e.\ $n=q_k=1$.  We are done
by induction on $k$.

This completes the proof of Lemma~\ref{lem:workhorse}.
\end{proof}

\subsection{Flow lines without trivial cylinders}
\label{sec:noCylinders}

We now put everything together to finish the proof of
Proposition~\ref{prop:braids}.  Let $C$ be a flow line without trivial
cylinders from $\{(\alpha_i,m_i)\}$ to $\{(\beta_j,n_j)\}$.  We start
with a local inequality for the incoming ends.  Recall that for each
$j$, the incoming ends of $C$ at $\beta_j$ determine a partition
$p_j^-=(q'_{j,1},q'_{j,2},\ldots)$ of $n_j$, and a braid $\xi_j^-$ in
a neighborhood of $\beta_j$ with $n_j$ strands.

\begin{lemma}[incoming inequality]
\label{lem:incomingInequality}
For a fixed $j$, write $\gamma=\beta_j$, $n=n_j$, $\xi=\xi_j^-$, and
$q_r=q'_{j,r}$.  Let $\tau\in\mc{T}(\gamma)$.  Then
\begin{equation}
\label{eqn:incomingInequality}
w_\tau(\xi)+\sum_r\mu_\tau(\gamma^{q_r})\ge\sum_{i=1}^n\mu_\tau(\gamma^i).
\end{equation}
Equality holds only if
$
\{q_1,q_2,\ldots\}=\pin(\gamma,n)
$.
\end{lemma}

\begin{proof}
The braid $\xi$ has components $\xi_1,\ldots,\xi_k$, where $\xi_i$ has
$q_i$ strands and winding number $\eta_i$.  Let
$\rho_i=\left\lceil\mu_\tau(\gamma^{q_i})/2\right\rceil$.  By
Lemma~\ref{lem:windingBound},
$
\eta_i\ge\rho_i
$.
Then
\begin{align}
\label{eqn:linkAndWrithe}
w_\tau(\xi) &= \sum_{i=1}^kw_\tau(\xi_i)+\sum_{i\neq
j}\ell_\tau(\xi_i,\xi_j)\\
\nonumber
&\ge \sum_{i=1}^k\rho_i(q_i-1)+\sum_{i\neq
j}\min(q_i\rho_j,q_j\rho_i).
\end{align}
The equality in the first line follows from the definitions of writhe
and linking number.  The second line follows from
Lemmas~\ref{lem:writheBound} and \ref{lem:linkingBound}.  Combining
this with Lemma~\ref{lem:workhorse}, we obtain the desired inequality
\eqref{eqn:incomingInequality}.

If equality holds in \eqref{eqn:incomingInequality}, then conditions
\ref{lem:writheBound}(i) and (ii) hold for each $q_i$, since we applied
Lemma~\ref{lem:writheBound} to each $\xi_i$, and conditions
\ref{lem:workhorse}(i) and (ii) hold.  All together, these conditions
imply that $\{q_1,q_2,\ldots\}=\pin(\gamma,n)$.
\end{proof}

By a symmetry argument, we can obtain a similar local inequality for
the outgoing ends.  If $C$ is a flow line as above, then for each
$i$, the outgoing ends of $C$ at $\alpha_i$ determine a partition
$p_i^+=(q_{i,1},q_{i,2},\ldots)$ of $m_i$, and a braid $\xi_i^+$ in
a neighborhood of $\alpha_i$ with $m_i$ strands.

\begin{lemma}[outgoing inequality]
\label{lem:outgoingInequality}
For a fixed $i$, write $\gamma=\alpha_i$, $m=m_i$, $\xi=\xi_i^+$, and
$q_r=q_{i,r}$.  Let $\tau\in\mc{T}(\gamma)$.  Then
\begin{equation}
\label{eqn:outgoingInequality}
w_\tau(\xi)+\sum_r\mu_\tau(\gamma^{q_r})\le\sum_{l=1}^m\mu_\tau(\gamma^l).
\end{equation}
Equality holds only if
$
\{q_1,q_2,\ldots\}=\pout(\gamma,m).
$
\end{lemma}

\begin{proof}
Let $Y'$ denote the mapping torus of $\phi^{-1}$.  By equation
\eqref{eqn:mappingTorus}, the map
\[
\begin{split}
\R\times (\Sigma\times\R)&\to \R\times (\Sigma\times\R),\\
(s,(x,t)) &\mapsto (-s,(x,-t))
\end{split}
\]
descends to a bijection
\begin{equation}
\label{eqn:symmetry}
\R\times Y \to \R\times Y'.
\end{equation}
This map sends the almost complex structure on $\R\times Y$ to an
admissible almost complex structure on $\R\times Y'$, see
\S\ref{sec:flowLines}.  For this almost complex structure on
$\R\times Y'$, the bijection \eqref{eqn:symmetry} sends flow lines to
flow lines.  Since the bijection \eqref{eqn:symmetry} switches $s$
with $-s$, incoming ends become outgoing ends and vice versa.  If
$\gamma$ is a periodic orbit in $Y$, let $\gamma'$ denote the
corresponding periodic orbit in $Y'$.  Then $C$ can be regarded as a
flow line in $\R\times Y'$ from $\{(\beta_j',n_j)\}$ to
$\{(\alpha_i',m_i)\}$.  By Lemma~\ref{lem:incomingInequality}, the
braid $\xi'$ at $\alpha_i'$ satisfies
\[
w_\tau(\xi')+\sum_r\mu_\tau({\gamma'}^{q_r})\ge
\sum_{l=1}^n\mu_\tau({\gamma'}^l),
\]
in the notation above equation \eqref{eqn:outgoingInequality}.
Because the bijection \eqref{eqn:symmetry} switches $t$ with $-t$, we
have $w_\tau(\xi')=-w_\tau(\xi)$ and
$\mu_\tau({\gamma'}^l)=-\mu_\tau(\gamma^l)$.  Thus we obtain the
inequality \eqref{eqn:outgoingInequality}.

If equality holds in \eqref{eqn:outgoingInequality}, then since we
used Lemma~\ref{lem:incomingInequality}, we must have
$\{q_1,\ldots,q_k\}=\pin(\gamma',m)$.  But
$\pin(\gamma',m)=\pout(\gamma,m)$; because if $\gamma$ is hyperbolic
with positive or negative eigenvalues, then so is $\gamma'$; and if
$\gamma$ is elliptic with monodromy angle $\theta$, then $\gamma'$ is
elliptic with monodromy angle $-\theta$.
\end{proof}

Proposition~\ref{prop:braids} follows immediately from
Lemmas~\ref{lem:incomingInequality} and \ref{lem:outgoingInequality}.
This proves the index theorem~\ref{thm:index} for flow lines without
trivial cylinders.

\section{Trivial cylinders}
\label{sec:cylinders}

To complete the proof of the index theorem \ref{thm:index}, we
need to consider a flow line $C$ which may contain trivial cylinders.
We begin with some notation, which will also be used in
\S\ref{sec:multiplyCovered}.  If $C$ is a flow line from $\alpha$ to
$\beta$, we write
$
\murel(C)=\murel(\alpha,\beta;[C]).
$
More generally, if $C_p$ is a flow line from the orbit set $\alpha_p$ to the
orbit set $\beta_p$, for
$p=1,\ldots,r$, we write
\[
\murel\left(\sum_{p=1}^rC_p\right)=
\murel\left(\sum_{p=1}^r\alpha_p,\sum_{p=1}^r\beta_p;\sum_{p=1}^r[C_p]\right).
\]
Here addition of orbit sets is defined by adding the multiplicities of
all periodic orbits involved.

Now to prove Theorem~\ref{thm:index}, write $C=C'\cup T$, where $T$ is
a union of trivial cylinders (possibly repeated), and $C'$ contains no
trivial cylinders, as in
\S\ref{sec:admissible}.  We have $\dim(\mc{M}_T)=0$; in fact
$\mc{M}_T=\{T\}$, by Proposition~\ref{prop:cylinderPoint}.  Therefore
\[
\dim(\mc{M}_C)=\dim(\mc{M}_{C'}).
\]
By Lemma~\ref{lem:indexBeginning} and Proposition~\ref{prop:braids},
\[
\dim(\mc{M}_{C'})\le I(C').
\]
So the following proposition will complete the proof of
Theorem~\ref{thm:index}.

\begin{proposition}[trivial cylinders]
\label{prop:cylinders}
If $C'$ is a flow line without trivial cylinders, and if $T$ is a
union of trivial cylinders (possibly repeated), then
\begin{equation}
\label{eqn:cylinderInequality}
\murel(C')\le\murel(C'+T)-2\#(C'\cap T).
\end{equation}
Equality holds only if $C=C'\cup T$ satisfies the admissibility
conditions \eqref{eqn:admissible1} and \eqref{eqn:admissible2}.
\end{proposition}

\begin{remark}
Of course, if $C$ is a flow line, then by definition $C'\cap
T=\emptyset$.  In general, by intersection positivity \cite{mcduff},
the algebraic intersection number $\#(C'\cap T)>0$ whenever $C'\cap
T\neq\emptyset$.
\end{remark}

\begin{lemma}[split partitions]
\label{lem:admissibleCylinders}
Let $\theta\in\R\setminus\Q$, and suppose
\begin{equation}
\label{eqn:fractionalPart}
\lfloor i\theta\rfloor + \lceil m\theta\rceil =  \lfloor (m+i)\theta\rfloor
\end{equation}
%(or equivalently $\{i\theta\}\ge \{(m+i)\theta\}$),
for $i=1,\ldots,n$. Then
\[
\pin(\theta,m+n)=\pin(m,\theta)\cup\pin(\theta,n).
\]
\end{lemma}

\begin{proof}
Write $\pin(\theta,m)=\{q_1,\ldots,q_k\}$ and
$\pin(\theta,n)=\{r_1,\ldots,r_l\}$.  By Lemma~\ref{lem:incoming},
we have $M_\theta(q_1,\ldots,q_k)=M_\theta(r_1,\ldots,r_l)=0$.  Using this, we
compute that
\[
\begin{split}
M_\theta(q_1,\ldots,q_k,r_1,\ldots,r_l) = &
2\sum_{i=1}^k\sum_{j=1}^l\min(q_i\lceil r_j\theta\rceil,r_j\lceil
q_i\theta\rceil) \\
&+2\sum_{i=1}^n(\lfloor i\theta\rfloor-\lfloor(m+i)\theta\rfloor).
\end{split}
\]
Using Lemma~\ref{lem:exercise}(a), we estimate
\[
\sum_{i=1}^k\sum_{j=1}^l\min(q_i\lceil r_j\theta\rceil,r_j\lceil
q_i\theta\rceil)
\le
\sum_{i=1}^k\sum_{j=1}^lr_j\lceil q_i\theta\rceil
= n\lceil m\theta\rceil.
\]
Putting this into the previous equation and then using
\eqref{eqn:fractionalPart}, we obtain
\[
M_\theta(q_1,\ldots,q_k,r_1,\ldots,r_l) \le 2\sum_{i=1}^n(\lfloor
i\theta\rfloor - \lfloor (m+i)\theta\rfloor + \lceil m\theta\rceil)
=0.
\]
By Lemma~\ref{lem:incoming}, it follows that
$\{q_1,\ldots,q_k,r_1,\ldots,r_l\}=\pin(\theta,m+n)$.
\end{proof}

\medskip
\noindent
{\em Proof of Proposition~\ref{prop:cylinders}.}  Suppose $T$ consists of
 trivial cylinders over periodic orbits $\gamma_l$ repeated $r_l$
 times.  At $\gamma_l$, suppose that $C'$ has outgoing ends of total
 multiplicity $m_l^+$, and incoming ends of total multiplicity
 $m_l^-$.  (Comparing with the notation of \S\ref{sec:admissible}, if
 $\gamma_l=\alpha_i=\beta_j$, then $m_l^+=m_i'$, $m_l^-=n_j'$, and
 $r_l=m_i-m_i'=n_j-n_j'$.)  Let $\xi_l^+$ and $\xi_l^-$ denote the
 corresponding braids from $C'$.  Let $\tau$ be a trivialization of
 $E$ over the ends of $C$.

Starting from $C$, we can perturb the trivial cylinders to obtain a
surface $S$ which is a representative of $[C]$ as in
Definition~\ref{def:rep}, except for finitely many self-intersections
coming from the intersections of $C'$ with $T$.  We have
\[
c_\tau(S)=c_\tau(C'),
\]
because the trivializations of $E$ are the same over both ends of each
cylinder.  By the relative adjunction formula \eqref{eqn:singular},
\[
\begin{split}
Q_\tau(S)-Q_\tau(C') - 2\#(C'\cap T)& = w_\tau(C')-w_\tau(S)\\ &=
2\sum_lr_l\left(\eta_\tau(\xi_l^-)-\eta_\tau(\xi_l^+)\right).
\end{split}
\]
And by definition of $\mu_\tau$, we have
\[
\mu_\tau(S)-\mu_\tau(C') = \sum_l\sum_{k=1}^{r_l}
\left(\mu_\tau\left(\gamma_l^{m_l^++k}\right) -
\mu_\tau\left(\gamma_l^{m_l^-+k}\right)\right).
\]
By the above three equations, in order to prove
Proposition~\ref{prop:cylinders}, it will suffice to show that for each $l$
and $1\le k\le r_l$,
\[
2(\eta_\tau(\xi_l^+)-\eta_\tau(\xi_l^-))
\le
\mu_\tau\left(\gamma_l^{m_l^++k}\right) -
\mu_\tau\left(\gamma_l^{m_l^-+k}\right),
\]
with equality for all such $k$ only if $C$ satisfies
\eqref{eqn:admissible1} and
\eqref{eqn:admissible2}.  In the rest of this proof, we drop the subscripts
`$l$' and `$\tau$'.

Let $(q_1^+,q_2^+,\ldots)$ and $(q_1^-,q_2^-,\ldots)$ denote the
partitions of $m^+$ and $m^-$ determined by the ends of $C'$ at
$\gamma$.  Let
$\rho_r^+=\left\lfloor\mu\left(\gamma^{q_r^+}\right)/2\right\rfloor$
and $\rho_s^- =
\left\lceil\mu\left(\gamma^{q_s^-}\right)/2\right\rceil$.
By Lemma~\ref{lem:windingBound}, we have $\eta(\xi^+)\le \sum_r\rho_r^+$ and
$\eta(\xi_i)\ge \sum_s\rho_s^-$.  (We know this even without assuming
$(\phi,J)$ is admissible at $\gamma$, by the asymptotics in
\cite{hwz2}.)  So it will suffice to show that for
$1\le k\le r$,
\begin{equation}
\label{eqn:cylinderGoal}
2\left(\sum_r\rho_r^+-\sum_s\rho_s^-\right)
\le
\mu\left(\gamma^{m^++k}\right) -
\mu\left(\gamma^{m^-+k}\right),
\end{equation}
with equality for all such $k$ only if \eqref{eqn:admissible1} and
\eqref{eqn:admissible2} hold for $C$ at $\gamma$.  To prove this, we
consider three cases.

\underline{Case 1}: $\gamma$ is hyperbolic with positive eigenvalues.
We can choose the trivialization $\tau$ so that $\mu(\gamma^k)=0$ for
all $k$.  Then all terms in \eqref{eqn:cylinderGoal} are zero, so
\eqref{eqn:cylinderGoal} holds.  The admissibility conditions
\eqref{eqn:admissible1} and \eqref{eqn:admissible2} are automatically
satisfied in this case.

\underline{Case 2}: $\gamma$ is hyperbolic with negative eigenvalues.  We can
choose the trivialization so that $\mu(\gamma^k)=k$ for all $k$.  The
inequality \eqref{eqn:cylinderGoal} asserts in this case that
\[
2\left(\sum_r\left\lfloor\frac{q_r^+}{2}\right\rfloor
-
\sum_s\left\lceil\frac{q_s^-}{2}\right\rfloor\right)
\le m^+-m^-.
\]
This clearly holds.  Equality holds only if $m^+$ and $m^-$ are even,
which is the content of
\eqref{eqn:admissible1} and
\eqref{eqn:admissible2} in this case, by
Remark~\ref{remark:admissible}, since $r>0$.

\underline{Case 3}: $\gamma$ is elliptic with monodromy angle $\theta$.  Then
\[
\sum_r\rho_r^+-\sum_s\rho_s^- = \sum_r\lfloor q_r^+\theta\rfloor -
\sum_s\lceil q_s^-\theta\rceil \le 
\lfloor m^+\theta\rfloor + \lfloor -m^-\theta\rfloor,
\]
and
\[
\mu\left(\gamma^{m^++k}\right) - \mu\left(\gamma^{m^-+k}\right)
= 2(\lfloor(m^++k)\theta\rfloor - \lfloor(m^-+k)\theta\rfloor).
\]
So in this case \eqref{eqn:cylinderGoal} follows from
\begin{equation}
\label{eqn:floors}
\lfloor m^+\theta\rfloor + \lfloor -m^-\theta\rfloor
\le \lfloor(m^+-m^-)\theta\rfloor
\le \lfloor (m^++k)\theta\rfloor - \lfloor(m^-+k)\theta\rfloor.
\end{equation}

Now suppose that equality holds in \eqref{eqn:cylinderGoal}, so that
it holds in \eqref{eqn:floors}, for $k=1,\ldots,r$.  Adding $\lfloor
k\theta\rfloor$ to both sides and manipulating, we obtain
\[
\lfloor m^+\theta\rfloor + \lfloor k\theta\rfloor -
\lfloor(m^++k)\theta\rfloor = \lceil m^-\theta\rceil + \lfloor
k\theta\rfloor - \lfloor(m^-+k)\theta\rfloor.
\]
The left side of this equation cannot be positive, and the right side
cannot be negative.  Hence both sides equal zero.  Since the right
side equals zero for $k=1,\ldots,n$, applying
Lemma~\ref{lem:admissibleCylinders} gives
\[
\pin(\theta,m^-+k) = \pin(\theta,m^-)\cup\pin(\theta,k).
\]
This is the admissibility condition \eqref{eqn:admissible2} for $C$ at
$\gamma$.  By symmetry, as in the proof of
Lemma~\ref{lem:outgoingInequality}, we also obtain the admissibility
condition \eqref{eqn:admissible1}.

This completes the proof of Proposition~\ref{prop:cylinders} and
the index theorem~\ref{thm:index}.
\qed

\section{Multiply covered pseudoholomorphic curves}
\label{sec:multiplyCovered}

Having proved the index theorem, we now want to establish compactness
results for moduli spaces of flow lines of dimension one and two.  A
major step is to show roughly that in these moduli spaces, a sequence
of flow lines cannot converge to a multiply covered pseudoholomorphic
curve, except that there can be repeated trivial cylinders.  We will
do that in \S\ref{sec:compactness}, by showing that an embedded (or
quasi-embedded) curve underlying such a multiply covered curve would
live in a moduli space of negative expected dimension, and hence does
not exist for generic $J$.  The key is the following index inequality;
see the beginning of \S\ref{sec:cylinders} for the notation.

\subsection{A generalization of the index inequality}

\begin{theorem}[multiply covered curves]
\label{thm:multiplyCovered}
Let $C_1,\ldots,C_r$ be disjoint flow lines and $d_1,\ldots,d_r$
positive integers.  Assume that $(\phi,J)$ is $d$-admissible, where
$d$ is the maximum period of the orbits at the ends of the $C_p$'s,
and $J$ is generic.  Then
\begin{equation}
\label{eqn:MCC}
\sum_{p=1}^rd_p\dim(\mc{M}_{C_p})
\le
\murel\left(\sum_{p=1}^rd_pC_p\right).
\end{equation}
\end{theorem}

\begin{remark}
Of course this theorem is a generalization of the index inequality
\eqref{eqn:indexInequality}, which is recovered when $r=d_1=1$.
We proved the latter separately in order to simplify the exposition.
\end{remark}

\begin{remark}
More generally, one might try to show that if $\alpha_p$ and $\beta_p$
are orbit sets for $p=1,\ldots,r$, and if $Z_p\in
H_2(Y;\alpha_p,\beta_p)$, then
\begin{equation}
\label{eqn:false}
\sum_{p=1}^rd_p\murel(\alpha_p,\beta_p;Z_p)
\le
\murel\left(\sum_{p=1}^rd_p\alpha_p,
\sum_{p=1}^rd_p\beta_p;\sum_{p=1}^rd_pZ_p\right).
\end{equation}
This, together with \eqref{eqn:indexInequality}, would imply
Theorem~\ref{thm:multiplyCovered}.  However, \eqref{eqn:false} is not
always true; one can use Proposition~\ref{prop:relativeIndex}(d) to
create a counterexample.
\end{remark}

\subsection{Topological preliminaries}

Before proving Theorem~\ref{thm:multiplyCovered}, we need some more
information about the relative intersection pairing.  If $C$ is a flow
line without trivial cylinders and $\tau$ is a trivialization of the
periodic orbits at the ends, define the total winding number
\[
\eta_\tau(C)=\sum_i\eta_\tau(\xi_i^+)-\sum_j\eta_\tau(\xi_j^-).
\]
Here $\xi_i^+$ and $\xi_j^-$ are the braids at the outgoing and
incoming ends, as in \S\ref{sec:adjunction}.

The following inequality will be needed in the proof of
Theorem~\ref{thm:multiplyCovered}, and is also useful for showing that
flow lines cannot exist under certain conditions.  For use in
\S\ref{sec:compactness}, instead of just considering flow lines, we
will more generally consider quasi-embedded GFL's, see
\S\ref{sec:GFL}.

\begin{proposition}[s-translation inequality]
\label{prop:s-trans}
If $C$ is a quasi-embedded GFL without trivial cylinders such that
$(\phi,J)$ is admissible near its ends, then
\[
Q_\tau(C)\ge -w_\tau(C)-\eta_\tau(C)+2\delta(C).
\]
\end{proposition}

\begin{proof}
Let $C'$ be the surface in $\R\times Y$ obtained by translating $C$ in
the $s$ direction, i.e.\ the $\R$ direction, by a small positive
amount.  Since the almost complex structure $J$ is $\R$-invariant,
$C'$ is pseudoholomorphic.  Let $C''=C\cup C'$.  Since $C$ contains no
trivial cylinders, $C''$ has no multiply covered components, so $C''$
is also a quasi-embedded GFL, see \S\ref{sec:GFL}.  We then have:
\[
\begin{split}
c_\tau(C'')&=2c_\tau(C),\\
Q_\tau(C'')&= 4Q_\tau(C),\\
\delta(C'')&=4\delta(C)+\#(C\cap C'),\\
w_\tau(C'')&=4w_\tau(C)+2\eta_\tau(C).
\end{split}
\]
The first three equations follow directly from the definitions.
To prove the fourth equation, we observe from the asymptotics in
\S\ref{sec:braids} that the braid from the outgoing (resp.\ incoming)
ends of $C'$ at $\gamma$ is obtained by moving the braid from $C$
outward from (resp.\ inward toward) $\gamma$.  Hence the braid for
$C''$ is the cabling of the braid from $C$ by a two-strand braid whose
writhe is twice the winding number.  Thus the fourth equation follows
from \eqref{eqn:writheOfCabling}.

We now apply equation \eqref{eqn:singular} to $C''$, subtract
equation \eqref{eqn:singular} applied to $C$ and $C'$,
and use the above four equations, to obtain
\[
Q_\tau(C)+w_\tau(C)+\eta_\tau(C)-2\delta(C)=\#(C\cap C').
\]
By intersection positivity \cite{mcduff}, $\#(C\cap C')\ge 0$.
\end{proof}

%For flow lines, there is a slightly simpler proof which I replaced.
% The essence is the following calculation:
%\[
%\begin{split}
%0 & \le \#(C\cap C')\\
%&= c_1(N,\pi_N\partial_s)\\
%&= c_1(N,\tau)+\eta_\tau(C)\\
%&=Q_\tau(C)+w_\tau(C)+\eta_\tau(C).
%\end{split}
%\]

Next, we introduce an extension of the relative intersection pairing
$Q$.  If $\alpha_p$ and $\beta_p$ are orbit sets for $p=1,2$, and if
$\tau$ is a trivialization of $E$ over all the periodic orbits in
$\alpha_1,\alpha_2,\beta_1,\beta_2$, then we can define
\[
Q_\tau:H_2(Y;\alpha_1,\beta_1)\times H_2(Y;\alpha_2,\beta_2)\to\Z
\]
just as in Definition~\ref{def:Q}.  Equivalently, if $Z_p\in
H_2(Y;\alpha_p,\beta_p)$ for $p=1,2$, then
\begin{equation}
\label{eqn:bilinear}
Q_\tau(Z_1+Z_2,Z_1+Z_2)=Q_\tau(Z_1,Z_1)+Q_\tau(Z_2,Z_2)+2Q_\tau(Z_1,Z_2).
\end{equation}

For this extended intersection pairing, there is a version of the
adjunction formula \eqref{eqn:adjunctionB}.  Let $S_1,S_2$ be
representatives of $Z_1,Z_2$ which intersect in finitely many points.
Let $\{\alpha_i\}$ (resp. $\{\beta_j\}$) denote the set of all
periodic orbits at which $S_1$ or $S_2$ has outgoing (resp.\ incoming)
ends.  For $p=1,2$, let $\xi_{p,i}^+$ (resp.\ $\xi_{p,j}^-$) denote
the braid corresponding to the outgoing (resp.\ incoming) ends of
$S_p$ at $\alpha_i$ (resp.\ $\beta_j$).  Define the total linking
number
\[
\ell_\tau(S_1,S_2)=\sum_i\ell_\tau(\xi_{1,i}^+,\xi_{2,i}^+)
- \sum_j\ell_\tau(\xi_{1,j}^-,\xi_{2,j}^-).
\]
This is symmetric in $S_1$ and $S_2$.  We then have:

\begin{lemma}
\label{lem:adjunctionC}
Under the above assumptions,
\[
Q_\tau(Z_1,Z_2)=-\ell_\tau(S_1,S_2)+\#(S_1\cap S_2).
\]
\end{lemma}

\begin{proof}
We can create $\tau$-representatives of $Z_1$ and $Z_2$, whose
projected conormals lie in distinct rays, by attaching cobordisms of
braids to $S_1$ and $S_2$ as in the proof of
Proposition~\ref{prop:adjunction}(b).  The intersection number of
these cobordisms equals minus the linking number.
\end{proof}

\subsection{Proof of Theorem~\ref{thm:multiplyCovered}}

We begin by assuming that the $C_p$'s do not contain trivial
cylinders.

Let $\{\alpha_i\}$ (resp.\ $\{\beta_j\}$) denote the set of all
periodic orbits at which any of the $C_p$'s have outgoing
(resp. incoming) ends.  Choose some trivialization $\tau$ of $E$ over
the $\alpha_i$'s and $\beta_j$'s.  Let $m_{p,i}$ (resp.\ $n_{p,j}$)
denote the total multiplicity of the outgoing (resp.\ incoming) ends of
$C_p$ at $\alpha_i$ (resp.\ $\beta_j$).  Let $M_i=\sum_pd_pm_{p,i}$
and $N_j=\sum_pd_pn_{p,j}$.  Write
$[C']=\sum_{p=1}^rd_p[C_p]$.  Extending the notation in
\S\ref{sec:begin}, write
\[
\begin{split}
c_\tau(C')&=c_\tau(E|_{[C']},\tau),\\
\mu_\tau(C')&=
\sum_i\sum_{k=1}^{M_i}\mu_{\tau_i^+}(\alpha_i^k)
-\sum_j\sum_{k=1}^{N_j}\mu_{\tau_j^-}(\beta_j^k),\\
Q_\tau(C')&=Q_\tau([C'],[C']).
\end{split}
\]
For the rest of this proof, we drop the subscript `$\tau$'.
By the definition of the relative index $\murel$,
\[
\murel(C')=c(C')+Q(C')+\mu(C').
\]
Now $c$ is linear, and $Q$ is bilinear:
\[
\begin{split}
c(C')&=\sum_{p=1}^rd_pc(C_p),\\
Q(C') &= \sum_{p=1}^rd_p^2Q(C_p) + \sum_{p\neq
p'}d_pd_{p'}Q(C_p,C_{p'}).
\end{split}
\]
By the index formula of Lemma~\ref{lem:indexBeginning} and our
genericity assumption on $J$,
\[
\dim(\mc{M}_{C_p})\le c(C_p)+Q(C_p)+w(C_p)+\mu^0(C_p).
\]
So by the above equations, our goal \eqref{eqn:MCC} is equivalent to
\begin{equation}
\label{eqn:MCC2}
\sum_{p=1}^r d_p(w(C_p)+\mu^0(C_p))
\le
\sum_{p=1}^r(d_p^2-d_p)Q(C_p) +
\sum_{p\neq p'} d_pd_{p'} Q(C_p,C_{p'})
+ \mu(C').
\end{equation}

We now eliminate $Q$ from this inequality.
By Proposition~\ref{prop:s-trans} and Lemma~\ref{lem:adjunctionC},
\[
\begin{split}
Q(C_p)& \ge -w(C_p)-\eta(C_p),\\
Q(C_p,C_{p'}) &=-\ell(C_p,C_{p'}).
\end{split}
\]
So to prove our goal \eqref{eqn:MCC2}, it suffices to show that
\begin{equation}
\label{eqn:MCC3}
\sum_{p=1}^r\left(d_p^2w(C_p) + (d_p^2-d_p)\eta(C_p) +
d_p\mu^0(C_p)\right) 
+ \sum_{p\neq p'}d_pd_{p'}\ell(C_p,C_{p'})
\le \mu(C').
\end{equation}

To prove \eqref{eqn:MCC3}, it will suffice to prove the following two
claims.  First, fix a periodic orbit $\gamma$, and suppose that the
flow line $C_p$ has incoming ends at $\gamma$ with multiplicities
$q_{p,1},\ldots,q_{p,k_p}$ totalling $n_p$.  Let $\xi_p$ denote the
resulting braid.  The first claim is that
\begin{equation}
\label{eqn:MCC4}
\begin{split}
\sum_{p=1}^r\left(d_p^2w(\xi_p) + (d_p^2-d_p)\eta(\xi_p) +
d_p\sum_{i=1}^{k_p}\mu(\gamma^{q_{p,i}})\right)&+\\
+ \sum_{p\neq p'}d_pd_{p'}\ell(\xi_p,\xi_{p'})
\ge&
\sum_{i=1}^{\sum_pd_pn_p}\mu(\gamma^i).
\end{split}
\end{equation}
Note that since we are now restricting attention to a single periodic
orbit, we are recycling the subscripts `$i$' here and `$j$' below.
The second claim is that an analogue of \eqref{eqn:MCC4} holds for
outgoing ends.  However by symmetry, as in the proof of
Lemma~\ref{lem:outgoingInequality}, it will suffice to prove only the
first claim \eqref{eqn:MCC4}.

To prove \eqref{eqn:MCC4}, define
\[
\rho_{p,i} = \left\lceil
\frac{\mu(\gamma^{q_{p,i}})}{2}\right\rceil.
\]
For $i=1,\ldots,k_p$, let $\eta_{p,i}$ denote the winding number with
respect to $\tau$ of the $i^{th}$ component of
$\xi_p$.  The winding bound of Lemma~\ref{lem:windingBound} gives
\begin{equation}
\label{eqn:etapi}
\eta_{p,i}\ge\rho_{p,i}.
\end{equation}
To simplify notation, let us now combine the two indices $p,i$ into a
single index $j$.  That is, we let each pair $(p,i)$ correspond to
some integer $j$, and we define $q_j=q_{p,i}$, $\rho_j=\rho_{p,i}$,
and $d_j=d_p$.  Using equation \eqref{eqn:linkAndWrithe},
Lemmas~\ref{lem:writheBound} and \ref{lem:linkingBound}, and the
inequality \eqref{eqn:etapi}, we obtain
\[
\begin{split}
\sum_p\left(d_p^2w(\xi_p) + (d_p^2-d_p)\eta(\xi_p)\right)
& + \sum_{p\neq p'}d_pd_{p'}\ell(\xi_p,\xi_{p'})\\
\ge &
\sum_{j,j'}d_jd_{j'}\min(q_j\rho_{j'},q_{j'}\rho_j) -
\sum_jd_j\rho_j.
\end{split}
\]
Thus our goal \eqref{eqn:MCC4} will follow from
\begin{equation}
\label{eqn:MCC5}
\sum_{j,j'}d_jd_{j'}\min(q_j\rho_{j'},q_{j'}\rho_j) +
\sum_j d_j\left(\mu(\gamma^{q_j})-\rho_j\right)
\ge
\sum_{i=1}^{\sum_jd_jq_j}\mu(\gamma^i).
\end{equation}
To prove \eqref{eqn:MCC5}, we simply apply the inequality of
Lemma~\ref{lem:workhorse} to the tuple
\[
(\underbrace{q_1,\ldots,q_1}_{\mbox{\scriptsize $d_1$ times}},\ldots,
\underbrace{q_j,\ldots,q_j}_{\mbox{\scriptsize $d_j$ times}},\ldots).
\]

This completes the proof of Theorem~\ref{thm:multiplyCovered} when the
$C_p$'s do not contain trivial cylinders.  We now have to show that if
$T$ is a union of trivial cylinders (possibly repeated), then
$\murel(C')\le \murel(C'+T)$.  The argument is a straightforward
modification of the proof of the inequality of
Proposition~\ref{prop:cylinders}.
\qed

\section{Compactness}
\label{sec:compactness}

We now prove the compactness theorem \ref{thm:compactness}.

\subsection{Rigidity of trivial cylinders}
\label{sec:rigidity}

We begin with a useful basic fact.

\begin{proposition}[rigidity]
\label{prop:cylinderPoint}
If $T$ is a union of trivial cylinders (possibly repeated),
then $T$ is the only flow line in its relative homology class (ending
at the same periodic orbits).
\end{proposition}

\begin{proof}
Let $u:C\to\R\times Y$ be a flow line homologous to $T$.  Since
$\int_T\omega=0$,
\[
\int_Cu^*\omega=0.
\]
By the definition of admissible almost complex structure, for any
tangent vector $v\in T(\R\times Y)$ we have
\begin{equation}
\label{eqn:omegaPositive}
\begin{split}
\omega(v,Jv)&\ge 0,\\
\omega(v,Jv)=0 &\Longleftrightarrow v\in\op{span}(\partial_s,\partial_t).
\end{split}
\end{equation}
Since $u$ is a pseudoholomorphic immersion by the definition of
flow line, it follows that it is locally a flat section of the
bundle $\R\times Y\to\R\times S^1$.  In other words, $u$ is a union of
trivial cylinders (possibly repeated), so $u=T$.
\end{proof}

\begin{remark}
A similar argument shows that in $S^1\times Y$, with an admissible
almost complex structure, any embedded pseudoholomorphic curve in an
$S^1$-invariant homology class is $S^1$ invariant.  This fact is
related to a result of Ionel and Parker \cite{ip}, and can be used
together with Taubes's ``SW=Gr'' theorem \cite{t1} to calculate the
Seiberg-Witten invariants of $Y$, see \cite{hl}.
\end{remark}

\subsection{The case $d=1$, $\op{genus}(\Sigma)>0$}

If $d=1$, then any flow line is a pseudoholomorphic section of the
bundle $\R\times Y\to\R\times S^1$.  The reason is that the fibers are
pseudoholomorphic and have intersection number 1 with the flow line,
and the flow line cannot contain a fiber or else it would not be
embedded; so it follows from intersection positivity \cite{mcduff}
that the flow line intersects each fiber transversely in a single
point.  If also $\op{genus}(\Sigma)>0$, then the compactness theorem
\ref{thm:compactness} is a special case of a standard result in Floer
theory; the assumption that $\op{genus}(\Sigma)>0$ implies that
$\pi_2(Y)=0$ so that there is no bubbling of pseudoholomorphic
spheres.

In \S\ref{sec:GFL}--\S\ref{sec:limits} we will prove
Theorem~\ref{thm:compactness} when $\partial\Sigma\neq\emptyset$ or
$d>\op{genus}(\Sigma)$.

\subsection{Generalized flow lines}
\label{sec:GFL}

To understand limits of flow lines, we need to consider more general
pseudoholomorphic curves which are not necessarily embedded.

\begin{definition}
\label{def:GFL}
A \define{generalized flow line} (GFL) from the orbit set
$\{(\alpha_i,m_i)\}$ to the orbit set $\{(\beta_j,n_j)\}$ is a
pseudoholomorphic map $u:C\to\R\times Y$, modulo reparametrization,
where:
\begin{itemize}
\item
$C$ is a punctured compact Riemann surface,
\item
$u$ has outgoing ends at $\alpha_i$ with total multiplicity $m_i$,
incoming ends at $\beta_j$ with total multiplicity $n_j$, and no other
ends.
\end{itemize}
\end{definition}
We write $\chi(u)=\chi(C)$ and $[u]=u_*[C]$.  We define the relative
index $\murel(u)=\murel(\alpha,\beta;[u])$.

We say that $u$ is a \define{quasi-embedding} if $u$ is an embedding
when restricted to the complement of a finite (possibly empty) set in
$C$.  If $u:C\to\R\times Y$ is any GFL and $(\phi,J)$ is admissible
near its ends, then $u$ factors through a quasi-embedding
$C'\to\R\times Y$ via a branched covering $C\to C'$, possibly the
identity covering.  The reason is that if $u$ has no multiply covered
components, then $u$ is a quasi-embedding, because the singularities
are isolated
\cite{mcduff}, and the discussion in \S\ref{sec:braids} shows that
there are no singularities for $|s|$ sufficiently large, where the
intersection of $u(C)$ with $\{s\}\times Y$ is an iterated nested
cabling of torus braids.

We say that $u$ is a \define{connector} if $u$ is a branched cover
(possibly the identity) of a union of trivial cylinders.  We say that
a GFL $u$ is \define{nontrivial} if $u$ is not a connector, or if $u$
is a connector with branch points.

Let $\widetilde{\mc{M}}_u$ denote the component of the moduli space of GFL's
containing $u$.  If $\tau$ is a trivialization of $E$ over the ends of
$u$, we define the \define{virtual dimension}
\begin{equation}
\label{eqn:vir-dim}
\op{vir-dim}(\widetilde{\mc{M}}_u)=2c_\tau(u)+\mu_\tau^0(u)-\chi(u).
\end{equation}
Here $\mu_\tau^0$ is defined as in \S\ref{sec:begin}.  The virtual
dimension is a homotopy invariant of $u$.  If $u$ is a quasi-embedding
we let $\widehat{\mc{M}}_u$ denote the component of the moduli space
of quasi-embedded GFL's containing $u$.  If $(\phi,J)$ is
$d$-admissible and $J$ is generic, if $u$ is a quasi-embedding, and if
$u$ does not contain any fibers of the projection $\R\times
Y\to\R\times S^1$ then $\widehat{\mc{M}}_u$ is cut out transversely by
Lemma~\ref{lem:transversality}(b).  In this case we further have
\begin{equation}
\label{eqn:vir}
\dim(\widehat{\mc{M}}_u)=\op{vir-dim}(\widetilde{\mc{M}_u}).
\end{equation}
This follows as in \cite{sft}, or by a modification of the proof of
Lemma~\ref{lem:indexBeginning}.

\begin{lemma}[parity of virtual dimension]
\label{lem:GFLParity}
If $\alpha$ and $\beta$ are admissible, and if $u$ is a GFL from
$\alpha$ to $\beta$, then
\[
\op{vir-dim}(\widetilde{\mc{M}}_u)\equiv I(u)\mod 2.
\]
\end{lemma}

\begin{proof}
We have
\[
I(u)-\op{vir-dim}(\widetilde{\mc{M}}_u) =
-c_\tau([u])+Q_\tau([u],[u])+\mu_\tau(u)-\mu_\tau^0(u)+\chi(u).
\]
Write $\alpha=\{(\alpha_i,m_i)\}$ and $\beta=\{(\beta_j,n_j)\}$.  By
equations \eqref{eqn:parity1} and \eqref{eqn:parity2},
\[
-c_\tau([u])+Q_\tau([u],[u])\equiv \sum_im_i-\sum_jn_j\mod 2.
\]
Let $k_i$ (resp.\ $l_j$) denote the number of ends of $u$ at
$\alpha_i$ (resp.\ $\beta_j$).  Then
\[
\chi(u)\equiv \sum_ik_i-\sum_jk_j\mod 2.
\]
Finally, one can see directly that since $\alpha$ and $\beta$ are
admissible,
\[
\mu_\tau(u)-\mu_\tau^0(u) \equiv \sum_i(m_i-k_i)-\sum_j(n_j-l_j)\mod
2.
\]
Combining the above four equations proves the lemma.
\end{proof}

\begin{lemma}[low index GFL's]
\label{lem:GFL}
Let $\alpha$ and $\beta$ be orbit sets of degree $d$, and assume that
$(\phi,J)$ is $d$-admissible and $J$ is generic.  Let $u:C\to\R\times
Y$ be a GFL from $\alpha$ to $\beta$.  Assume that
$\partial\Sigma\neq\emptyset$ or $d>\op{genus}(\Sigma)$.  Then:
\begin{description}
\item{(a)}
$\murel(u)\ge 0$.
\item{(b)}
If $\murel(u)=0$, then $u$ is a connector.
\item{(c)}
If $\murel(u)=1$, or $\murel(u)=2$ and $\alpha,\beta$ are admissible,
then $u$ is the disjoint union of a nontrivial flow line and a
(possibly empty) connector.
\end{description}
\end{lemma}

\begin{proof}
Let $C_1,\ldots,C_r$ be the components of the quasi-embedding $C'$
above, regarded as subsets of $\R\times Y$, and let $d_p$ denote the
covering multiplicity of $C$ over $C_p$.  The projection $\R\times
Y\to\R\times S^1$ restricts to a holomorphic map on each component
$C_p$.  Therefore:
\begin{description}
\item{(*)}
Each $C_p$ either has at least one
incoming and one outgoing end, or is a fiber of the projection
$\R\times Y\to\R\times S^1$.
\end{description}
We now consider two cases.

$\underline{\mbox{Case 1:}}$ Suppose that no $C_p$ is a fiber.  Then
\begin{equation}
\label{eqn:mcg}
\sum_{p=1}^rd_p\dim(\widehat{\mc{M}}_{C_p})\le \murel(u)-2\Delta,
\end{equation}
where $\Delta$ is a nonnegative integer which is zero only if the
$C_p$'s are embedded and disjoint.  In fact, \eqref{eqn:mcg} holds
with
\[
\Delta=\sum_pd_p^2\delta(C_p)+\sum_{p<p'}d_pd_{p'}\#(C_p\cap C_{p'}).
\]
The proof of \eqref{eqn:mcg} follows the proof of
Theorem~\ref{thm:multiplyCovered}, except that we use equations
\eqref{eqn:vir}, \eqref{eqn:vir-dim}, and \eqref{eqn:singular} in
place of Lemma~\ref{lem:indexBeginning}.

Assertion (a) now follows from \eqref{eqn:mcg}.  To prove (b), if
$\murel(u)=0$, then $\dim(\widehat{\mc{M}}_{C_p})=0$ for each $p$ by
\eqref{eqn:mcg}, so $C_p$ is fixed under the translation action of
$\R$ on $\R\times Y$, so $C_p$ is a union of trivial cylinders.

To prove (c), we first observe that if $\murel(u)\in\{1,2\}$, then the
$C_p$'s are embedded and disjoint, or else by $\eqref{eqn:mcg}$ we
would have $\dim(\widehat{\mc{M}}_{C_p})=0$ so that the $C_p$'s
would be unions of trivial cylinders, giving $\murel(u)=0$.

Now suppose that $\murel(u)=1$.  Then for some $p$, we have
$\dim(\widehat{\mc{M}}_{C_p})=d_p=1$.  So $C_p$ is a nontrivial flow
line; the other components of $u$ are connectors as in part (b), and
they are disjoint from $C_p$ as explained above.

Finally suppose that $\murel(u)=2$ and $\alpha,\beta$ are admissible.
We claim that if $C_p$ is nontrivial then $d_p=1$, so we are done as
in the case $\murel(u)=1$.  To prove the claim, suppose $d_p>1$.
Since $\alpha$ and $\beta$ are admissible, all ends of $C_p$ are
elliptic.  It follows that $\dim(\widehat{\mc{M}}_{C_p})$ is even, by
equations \eqref{eqn:vir-dim} and \eqref{eqn:vir},
Lemma~\ref{lem:GFLParity}, and
Proposition~\ref{prop:relativeIndex}(c).  By \eqref{eqn:mcg},
$\dim(\widehat{\mc{M}}_{C_p})\le 1$, so we in fact have
$\dim(\widehat{\mc{M}}_{C_p})= 0$.  Therefore $C_p$ is a union of
trivial cylinders.

$\underline{\mbox{Case 2:}}$ Suppose that some $C_p$ is a fiber of the
projection $\R\times Y\to\R\times S^1$.  Then
$\partial\Sigma=\emptyset$, so by assumption
$d>g=\op{genus}(\Sigma)$.  Let $u'$ denote the GFL obtained from $u$
by deleting the component(s) covering $C_p$.  By the index ambiguity
formula of Proposition~\ref{prop:relativeIndex}(d), we have
\[
I(u')=I(u)-2d_p(d-g+1)\le I(u)-4.
\]
Thus removing all fibers from $u$ decreases the relative index by at
least 4.  So by assertion (a) in Case 1, $I(u)\ge 4$.  Thus assertion
(a) holds in this case as well, while assertions (b) and (c) are
vacuously true.
\end{proof}

\subsection{Limits of sequences of flow lines}
\label{sec:limits}

To describe the possible limits of sequences of flow lines, we need
the following definitions.

\begin{definition}
A \define{$k$-times broken GFL} from $\alpha$ to $\beta$ is a sequence
$(u_1,\ldots,u_k)$, where there exist orbit sets $\alpha_0,\ldots,
\alpha_k$ such that $u_i$ is a nontrivial GFL
from $\alpha_{i-1}$ to $\alpha_i$, and $\alpha_0=\alpha$ and
$\alpha_k=\beta$.  We also assume that at each periodic orbit in
$\alpha_i$, the multiplicities of the incoming ends of $u_i$ agree
with the multiplicities of the outgoing ends of $u_{i+1}$.
\end{definition}

If $\psi\in\R$, let $T_\psi:\R\times Y\to\R\times Y$ denote the translation
sending $(s,y)\mapsto (s-\psi,y)$.

\begin{definition}
\label{def:convergence}
Let $(u_1,\ldots,u_k)$ be a broken GFL from $\alpha$ to $\beta$.  Let
$\{v_1,v_2,\ldots\}$ be a sequence of GFL's from $\alpha$ to $\beta$
in a fixed relative homology class $Z$ whose domains have a fixed
topological type.  We say that $\{v_n\}$ {\em converges\/} to
$(u_1,\ldots,u_k)$ if:
\begin{description}
\item{(a)}
The domain of $v_n$ converges in Deligne-Mumford space (possibly after
adding some marked points to the $v_n$'s and $u_i$'s) to a nodal curve
which can be decomposed into a union of curves $C_1,\ldots,C_k$, such
that $C_i$ intersects $C_{i+1}$ in one node for each incoming end of
$u_i$ and $C_{i-1}$ in one node for each outgoing end of $u_i$; and
after puncturing $C_i$ at these points, we obtain the domain of $u_i$,
possibly with some pairs of points identified to nodes.
\item{(b)}
There are real numbers
$s_{1,n}>s_{2,n}>\cdots>s_{k,n}$ for each positive integer $n$ such
that $T_{s_{i,n}}\circ v_n\to u_i$ in $C^\infty$ on compact
sets.
\item{(c)}
$\sum_{i=1}^k[u_i]=Z$.
\end{description}
\end{definition}

We now have the following version of Gromov compactness.  If
$u:C\to\R\times Y$ is a pseudoholomorphic map, we define the {\em
energy} $E(u)=\int_Cu^*\omega$.

\begin{lemma}[general compactness]
\label{lem:gc}
Let $\{C_n\}_{n=1,2,\ldots}$ be a sequence of flow lines from $\alpha$
to $\beta$, with the energy $E(C_n)$ uniformly bounded from above.
Assume that $(\phi,J)$ is admissible near the periodic orbits in
$\alpha$ and $\beta$.  Then there is a subsequence with a fixed
relative homology class and topological type which converges in the
above sense to a broken GFL $(u_1,\ldots,u_k)$.
\end{lemma}

The proof of Lemma~\ref{lem:gc} is mostly a standard argument using
Gromov compactness, except that we need to be careful because the
genus of $C_n$ is not a priori bounded.  Fortunately, in dimension
four a version of Gromov compactness is available using currents which
does not assume any such bound \cite{t98,y}.  We will use the
following special case of it.

\begin{lemma}[currents]
\label{lem:taubes}
Fix $d\in\Z$ and $E_0\in\R$.  Let $u_k:C_k\to\R\times Y$ be a sequence
of proper pseudoholomorphic maps such that for each $k$, the
current $u_k[C_k]$ is a homology between two orbit sets of degree
$d$, and $E(u_k)<E_0$.  Then we can pass to a subsequence such that:
\begin{itemize}
\item
The $u_k$'s converge weakly as currents in $\R\times Y$ to a
proper pseudoholomorphic map $u:C\to\R\times Y$.
\item
For any compact set $K\subset\R\times Y$, as $k\to\infty$,
\begin{equation}
\label{eqn:geometricConvergence}
\sup_{x\in u_k(C_k)\cap K}\op{dist}(x,u(C))+\sup_{x\in
u(C)\cap K}\op{dist}(x,u_k(C_k))\to 0.
\end{equation}
\end{itemize}
\end{lemma}

\begin{proof}
We first note that there is a natural symplectic form
\[
\Omega=\omega+ds\wedge dt
\]
on $\R\times Y$.  An admissible almost complex structure $J$ on $Y$ is
tamed by $\Omega$:
\begin{equation}
\label{eqn:tame}
\mbox{$\Omega(v,Jv)>0$ for $v\neq 0$.}
\end{equation}
Furthermore, $\R\times Y$ is exhausted by compact sets $[a,b]\times
Y$, on each of which we have a uniform bound
\[
\int_{u_k[C_k]\cap ([a,b]\times Y)}\Omega \le d(b-a)+E_0
\]
by \eqref{eqn:omegaPositive}.  The lemma now follows as in
\cite[Prop.\ 3.3]{t98} or \cite{y}.
\end{proof}

\medskip

\noindent
{\em Proof of Lemma~\ref{lem:gc}.}  We proceed in 4 steps.  As usual
let $d$ denote the degree of $\alpha$ and $\beta$.  We choose a metric
on $Y$, and work with the product metric on $\R\times Y$.  We can pass
to a subsequence such that $E(C_n)\to E_0$.  Let $\Gamma_d$ denote the
union of the periodic orbits of period $\le d$.

\underline{Step 1}:
We claim that for any sequence of real numbers $s_n$, we can pass to a
subsequence such that $T_{s_n}(C_n)$ converges in the sense of
Lemma~\ref{lem:taubes} to a GFL.

To prove this, we can pass to a subsequence such that $T_{s_n}(C_n)$
converges in the sense of Lemma~\ref{lem:taubes} to some proper
pseudoholomorphic map $u:C\to\R\times Y$.  Since the integral of
$u^*\omega$ is locally nonnegative, using the weak convergence of
currents we obtain $E(u)\le E_0$.

Now $u(C)$ is $C^0$-asymptotic to a union of trivial cylinders as
$s\to+\infty$.  Otherwise there exists $\epsilon>0$ and a sequence
$\psi_k\to\infty$ such that $u(C)\cap(\{\psi_k\}\times Y)$ contains a
point of distance greater than $\epsilon$ from $\R\times\Gamma_d$.
The sequence of translates $T_{\psi_k}\circ u$ contains a subsequence
converging in the sense of Lemma~\ref{lem:taubes} to some
pseudoholomorphic curve.  Since $u$ has finite energy and
$\psi_k\to\infty$, the limiting curve has energy zero and hence its
image is a union of trivial cylinders as in
Proposition~\ref{prop:cylinderPoint}.  Together with
\eqref{eqn:geometricConvergence}, this contradicts our
distance assumption.  Likewise, $u(C)$ is $C^0$-asymptotic to a
possibly different union of trivial cylinders as $s\to -\infty$.

We will see below, using the assumption that the $C_n$'s are flow
lines, that we can pass to a further subsequence so that the
topological type of $C_n$ is fixed.  Then the domain $C$ is a
punctured compact Riemann surface, and by standard lemmas we have
$C^\infty$ convergence to trivial cylinders at the punctures, so $u$
is a GFL.

\underline{Step 2}:  Suppose $u:C\to\R\times Y$ is a proper
pseudoholomorphic map such that $u[C]$ is a homology between orbit
sets of degree $d$, and suppose $u(C)$ is not a union of trivial
cylinders.  Then we claim that there is a constant $\delta>0$
depending only on $d$ such that $E(u)>\delta$.

If not, then we can take a sequence $v_k:C'_k\to\R\times Y$ of such
maps with $E(v_k)\to 0$.  Let $N$ denote the set of points within
distance $\epsilon$ of $\Gamma_d$, where $\epsilon>0$ is chosen small
enough so that $N$ is a tubular neighborhood of $\Gamma_d$.  Since
$v_k(C'_k)$ is not a union of trivial cylinders, $E(v_k)>0$ by
equation~\eqref{eqn:omegaPositive}, so on homological grounds
$v_k(C'_k)\not\subset\R\times N$.  By applying translations we may
assume that $v_k(C'_k)\cap(\{0\}\times Y)$ contains a point of distance
at least $\epsilon$ from $\R\times \Gamma_d$.  Applying
Lemma~\ref{lem:taubes} as in the third paragraph of Step 1 gives a
contradiction.

\underline{Step 3}:
We now apply Step 1 with judicious choices of $s_n$.  We can assume
that $C_n$ is not a union of trivial cylinders for large $n$, or else
the lemma is trivially true.  We can then define
\[
s_{1,n}=\sup\{s\in\R\mid C_n\cap(\{s\}\times Y)\not\subset
\{s\}\times N\}.
\]
By Step 1, we can pass to a subsequence so that $T_{s_{1,n}}(C_n)$
converges to a GFL $u_1$, whose image is not a union of trivial
cylinders, with $E(u_1)\le E_0$.

Suppose $E(u_1)<E_0$.  Since $u_1$ is asymptotic to some union of
periodic orbits of period $\le d$ as $s\to-\infty$, we can find
$\psi_1\in\R$ such that $\op{dist}(y,\Gamma_d)<\epsilon/2$ whenever
$y\in Y$, $(s,y)\in\op{Im}(u_1)$, and $s<\psi_1$.  Then for large $n$,
$C_n\cap ((-\infty,s_{1,n}+\psi_1)\times Y)$ is not contained in
$\R\times N$, or else $C_n$ would be homologous to $u_1$ and would
have the same energy.  We can then define
\[
s_{2,n}=\sup\{s<s_{1,n}+\psi_1\mid C_n\cap (\{s\}\times Y) \not\subset
\{s\}\times N\}.
\]
By Step 1, we can pass to a subsequence so that $T_{s_{2,n}}(C_n)$
converges to a GFL $u_2$.

Continuing this process, suppose that $s_{i,n}$ have been chosen for
$i<k$.  If $E(u_1)+\cdots+E(u_{k-1})<E_0$, we define $\psi_{k-1}$ from
$u_{k-1}$ as above, set
\[
s_{k,n}=\sup\{s<s_{k-1,n}+\psi_{k-1}\mid C_n\cap (\{s\}\times Y) \not\subset
\{s\}\times N\},
\]
and pass to a subsequence so that $T_{s_{k,n}}(C_n)$ converges to a
GFL $u_k$.

This process must eventually stop since the sum of the $E(u_i)$'s is
bounded from above by $E_0$, while by Step 2 each $E(u_i)$ is bounded
from below by $\delta>0$.  We obtain GFL's $u_1,\ldots,u_k$ and real
numbers $s_{1,n}>s_{2,n}>\cdots>s_{k,n}$ for each $n$ such that the
relative homology class $[C_n]$ is eventually constant and equals
$\sum_i[u_i]$, and $T_{s_{i,n}}(C_n)\to u_i$ in the sense of
Lemma~\ref{lem:taubes}.

\underline{Step 4}:
We can pass to a subsequence so that the relative homology class
$[C_n]$ is fixed.  By Theorem~\ref{thm:euler}(a), there is a lower
bound on the Euler characteristic $\chi(C_n)$.  This implies that we
can pass to a subsequence so that the topological type of $C_n$ is
fixed.  As in \cite{ehwz}, we can then obtain convergence in the sense
of Definition~\ref{def:convergence}, by passing to a further
subsequence and possibly inserting some nontrivial connectors in
between $u_1,\ldots,u_k$ to obtain the limiting broken GFL.
\qed

Together with Lemma~\ref{lem:gc}, the following lemma will complete
the proof of Theorem~\ref{thm:compactness}(a).

\begin{lemma}[possible limits]
\label{lem:broken}
Let $\alpha$ and $\beta$ be orbit sets of degree $d$.  Let
$(u_1,\ldots,u_k)$ be a broken GFL from $\alpha$ to $\beta$ which is
the limit of a sequence of flow lines $\{C_n\}$ in the same relative
homology class of relative index $I_0\in\{1,2\}$.  Assume that
$\partial\Sigma\neq\emptyset$ or $d>\op{genus}(\Sigma)$.  Suppose
$(\phi,J)$ is $d$-admissible and $J$ is generic.  Then:
\begin{description}
\item{(a)}
If $I_0=1$, then $k=1$ and $u_1$ is a nontrivial flow line.
\item{(b)}
If $I_0=2$ and $\alpha,\beta$ are admissible, then $u_1$ and $u_k$ are
nontrivial flow lines, and $u_i$ is a nontrivial connector for
$1<i<k$.
\end{description}
\end{lemma}

\begin{proof}
Assume that (a) $I_0=1$, or (b) $I_0=2$ and
$\alpha,\beta$ are admissible.

\underline{Step 1}: By the additivity of the relative index
in Proposition~\ref{prop:relativeIndex}(b), we have
$
\sum_{i=1}^kI(u_i)=I_0.
$
By Lemma~\ref{lem:GFL}, it follows that each $u_i$ is either a
nontrivial connector, or the disjoint union of a nontrivial flow line
and a (possibly empty) connector.
%(Note that if $I(u_i)=2$, then the
%ends of $u_i$ are admissible, so Lemma~\ref{lem:GFL}(c) does apply;
%otherwise $u_j$ is not a connector for some $j\neq i$, and then
%$I_0>2$.)

\underline{Step 2}: 
The GFL's $u_1,\ldots,u_k$ cannot all be connectors, or else $I(C_n)=0$.

\underline{Step 3}:
A connector lives in a moduli space of nonnegative virtual dimension,
see Step 5.  A nontrivial flow line lives in a moduli space of
positive virtual dimension, by equation \eqref{eqn:vir} and the
$\R$-action.  Also, by \eqref{eqn:vir-dim}, the virtual dimension is
additive:
\[
\sum_{i=1}^k\op{vir-dim}(\widetilde{\mc{M}}_{u_i})
=\lim_{n\to\infty}\op{vir-dim}(\widetilde{\mc{M}}_{C_n}) \le I_0.
\]
The inequality on the right holds by the relative adjunction formulas
and the inequality \eqref{eqn:braids}.

\underline{Step 4}: The $C_n$'s are admissible.
Otherwise, if some $C_n$ is not admissible, then
$\op{vir-dim}(\widetilde{\mc{M}}_{C_n})<I_0$ as in
Theorem~\ref{thm:index}, and moreover in case (b) the difference is at
least two by Lemma~\ref{lem:GFLParity}.  Consequently
$\dim(\widehat{\mc{M}}_{C_n})=0$ by \eqref{eqn:vir}, so $C_n$ is a
union of trivial cylinders, so $I(C_n)=0$, a contradiction.

\underline{Step 5}: $u_1$ and $u_k$ do not contain any nontrivial connectors.

Otherwise, WLOG $u_k$ contains a nontrivial
connector $v$.  Then $\op{vir-dim}(\widetilde{\mc{M}}_v)\ge 1$, and in
case (b) $\op{vir-dim}(\widetilde{\mc{M}}_v)\ge 2$.  The reason is
that we can compute $\op{vir-dim}(\widetilde{\mc{M}}_v)$ from equation
\eqref{eqn:vir-dim} by starting at the bottom of $v$ and moving
vertically past one branch point at a time.  By homotopy invariance of
virtual dimension, we can assume that the branch points are simple.
Then each time we move past a branch point, two circles in $v$ join
into one or one circle separates into two.  By
Proposition~\ref{prop:CZ}, the contribution to the virtual dimension
is 1 for hyperbolic orbits, and 0 or 2 for elliptic orbits.  If the
branching is over an elliptic orbit, then the bottom branch point
contributes 2 to the virtual dimension, by Step 4 and
Lemma~\ref{lem:exercise}(b),(c).  In case (b), the branching is always
over an elliptic orbit, since $\alpha$ and $\beta$ are admissible.

By Step 3, all other components of the $u_i$'s live in moduli spaces
of virtual dimension $0$, so no component of any $u_i$ is a nontrivial
flow line.  But this contradicts Steps 1 and 2.

\underline{Step 6}: Now we are done.  By Steps 1 and 5, $u_1$ and
$u_k$ are nontrivial flow lines.  By Step 3, $u_i$ is a nontrivial
connector for $1<i<k$.  If $I_0=1$, then we must have $k=1$ by
additivity of the relative index.
\end{proof}

\subsection{Fiber bubbles and transversality}
\label{sec:almostAdmissible}

The proof of Theorem~\ref{thm:compactness}(a) is now complete.  We now
prove Theorem~\ref{thm:compactness}(b).  By
Theorem~\ref{thm:compactness}(a), we can assume that
$\partial\Sigma=\emptyset$ and $1<d\le g$.  In this case, compactness
might not hold for an admissible almost complex structure, because the
fibers of the projection $\R\times Y\to\R\times S^1$ are
pseudoholomorphic and might bubble off.  To ensure compactness, we need
to relax the requirement that the almost complex structure preserve
$E$.

\begin{definition}
\label{def:almost}
Let $J$ be an almost complex structure on $\R\times Y$.  We say that
$(\phi,J)$ is {\em almost $d$-admissible\/} if:
\begin{description}
\item{(a)}
$J$ is $\R$-invariant.
\item{(b)}
$J(\partial_s)=\partial_t$.
\item{(c)}
There is an almost complex structure $J'$ such that $(\phi,J')$ is
$d$-admissible and $J$ agrees with $J'$ on $\R$ cross a neighborhood
of the periodic orbits of period $\le d$ and a neighborhood of
$\partial Y$.
\item{(d)}
$J$ is $\Omega$-tame, i.e.\ equation \eqref{eqn:tame} holds.
\end{description}
\end{definition}
Of course, if $J$ is $d$-admissible, then $J$ is almost
$d$-admissible.  Also note that if $J$ is almost $d$-admissible, then
a $C^0$-small perturbation (with respect to an $\R$-invariant metric
on $\R\times Y$) of $J$ satisfying conditions (a), (b), and (c) will
automatically satisfy the taming condition (d).

Condition (b) implies that $J$ preserves a 2-plane bundle
$E'\subset TY$, and (b) and (d) imply that $E'$ is homotopic to $E$.
This implies that the relative adjunction formula
\eqref{eqn:adjunctionA} still holds.

Condition (b) implies that \eqref{eqn:omegaPositive} still
holds, so Proposition~\ref{prop:cylinderPoint} still holds.

\begin{lemma}[transversality]
\label{lem:transversality}
Let $C$ be a quasi-embedded GFL and let $d$ be any positive integer.
Then $\widehat{\mc{M}}_C$ is cut out transversely, provided that
either:
\begin{description}
\item{(a)}
$(\phi,J)$ is almost $d$-admissible and $J$ is generic, or
\item{(b)}
$(\phi,J)$ is $d$-admissible and $J$ is generic, and $C$ contains no
fibers of the projection $\R\times Y\to\R\times S^1$.
\end{description}
\end{lemma}

\begin{proof}
(a) We can assume that our GFL's contain no trivial cylinders, as
these are always cut out transversely.  The proof is now a slight
modification of the proof of \cite[Thm.\ 5.1(i)]{fhs}.
%Recall that for an immersion $u:C\to\R\times Y$, the
%pseudoholomorphic curve equation can be regarded as a well-defined
%section
%\[
%\overline{\partial}_J(u)\in\Hom^{0,1}_J(TC,N_u),
%\]
%which vanishes if and only if $u$ is pseudoholomorphic.
The essential point is the following:

Claim.  For any almost $d$-admissible $(\phi,J)$ and any
quasi-embedded GFL $C$ without trivial cylinders, there is a nonempty
open set $U\subset C$ away from a neighborhood of the period $\le d$
periodic orbits, such that for each $x\in U$:
\begin{description}
\item{(i)}
$\pi^{-1}(\pi(x))=\{x\}$, where $\pi:C\to Y$ denotes the projection.
\item{(ii)}
$C$ is nonsingular at $x$, and the derivative
of $\overline{\partial}_J(C)$ with respect to $J$,
namely the bundle map
\begin{equation}
\label{eqn:restrictionMap}
\left\{\psi\in\Hom^{0,1}_J(TX|_C,TX|_C)\Big|
\psi|_{\op{span}(\partial_s,\partial_t)}=0\right\}
\longrightarrow
\Hom^{0,1}_J(TC,N_C),
\end{equation}
defined by restricting $\psi$ to $TC$ and projecting to $N_C$, is
surjective at $x$.  Here $\Hom^{0,1}_J$ denotes the space of
$J$-antilinear maps.
\end{description}
%  Note
%that condition (i) is needed because of the $\R$-invariance condition
%(a) above, and the constraint on the domain of the map
%\eqref{eqn:restrictionMap} arises from condition (b) above.)

To prove the claim, we will in fact obtain properties (i) and (ii) on
an open dense set $U$ in $C$.  Condition (i) holds on a dense open set
just as in \cite[Thm.\ 1.13]{hwz3}.  Regarding (ii), we observe that
with respect to $J$, the bundle map \eqref{eqn:restrictionMap} is
complex linear.  Since the bundle $\Hom^{0,1}_J(TC,N_C)$ has complex
rank one, the map \eqref{eqn:restrictionMap} is surjective except
where it is zero.  But it is zero only where $C$ is tangent to
$\op{span}(\partial_s,\partial_t)$.  The set of points where $C$ is
not tangent to $\op{span}(\partial_s,\partial_t)$ is certainly open,
and it is also dense, or else by a unique continuation argument $C$
would contain a trivial cylinder.

(b) We first note that if $C$ does not contain a fiber, then no GFL in
$\widehat{\mc{M}}_C$ contains a fiber, because any pseudoholomorphic
curve close to a fiber is also a fiber, as the projection $\R\times
Y\to\R\times S^1$ is holomorphic.  The proof of transversality is now
similar to part (a).  In this case we assume that $C$ contains no
trivial cylinders and fibers; and we show that on an open dense set
$U\subset C$, for $x\in U$, conditions (i) and (ii) above hold, but
with the map
\eqref{eqn:restrictionMap} replaced by the map
\begin{equation}
\label{eqn:restrictionMap2}
\Hom^{0,1}_J(E|_C,E|_C)
\longrightarrow
\Hom^{0,1}_J(TC,N_C)
\end{equation}
that sends
$\psi\in\Hom^{0,1}_J(E|_C,E|_C)$ to the composition
\[
TC
\longrightarrow
E|_C
\stackrel{\psi}{\longrightarrow}
E|_C
\longrightarrow
N|_C.
\]
Here the left and right arrows are projections.  The map
\eqref{eqn:restrictionMap2} is surjective at $x$ unless $T_xC$ agrees
with $\op{span}(\partial_s,\partial_t)$ or $E_x$.  Thus, by a unique
continuation argument, the map \eqref{eqn:restrictionMap2} is
surjective on an open dense set in $C$ as long as $C$ does not contain
a fiber or a trivial cylinder.
\end{proof}

\begin{lemma}
If $(\phi,J)$ is almost $d$-admissible and $J$ is generic, then
Theorem~\ref{thm:index} still holds, and the inequality
\eqref{eqn:mcg} holds for any quasi-embedded GFL's $C_1,\ldots,C_r$.
\end{lemma}

\begin{proof}
All necessary changes to the proofs are handled by
Lemma~\ref{lem:transversality}(a) and the remarks preceding it.
The rest of the proofs of Theorem~\ref{thm:index} and the inequality
\eqref{eqn:mcg} either do not involve $J$, or are local to the
periodic orbits of period $\le d$ and hence are unchanged by condition
(c) in Definition~\ref{def:almost}.
\end{proof}

\begin{lemma}
\label{lem:last}
There exists an almost $d$-admissible $(\phi,J)$ such that $J$ is
generic as in Lemma~\ref{lem:transversality}(a), for which the
analogue of Lemma~\ref{lem:GFL} holds, assuming only $g>1$ instead of
$\partial\Sigma\neq\emptyset$ or $d>\op{genus}(\Sigma)$.
\end{lemma}

\begin{proof}
We have to modify the proof of Lemma~\ref{lem:GFL} in two ways.

First, statement (*) in that proof is not necessarily true, as the
projection $\R\times Y\to\R\times S^1$ may no longer be
pseudoholomorphic.  Nonetheless, there exists almost $d$-admissible
$(\phi,J)$ with $J$ generic such that instead we have:
\begin{description}
\item{(*$'$)}
Each $C_p$ either has at least one incoming and one outgoing end, or
is homologous to a positive multiple of the class of the fiber of the
projection $\R\times Y\to\R\times S^1$.
\end{description}
If not, then we can take a sequence of generic $J_n$ with $(\phi,J_n)$
almost $d$-admissible, converging to a generic $J_\infty$ such that
$(\phi,J_\infty)$ is $d$-{\em admissible\/}, together with
$J_n$-pseudoholomorphic curves that do not satisfy (*$'$).  By Gromov
compactness, a subsequence of these curves converges to a
$J_\infty$-holomorphic curve which by (*) cannot exist.

Second, we have to show that $C_p$ cannot to be homologous to $k$
times the fiber class.  If it is, then using equation
\eqref{eqn:singular} and our assumption that $g>1$, we have
\[
\op{vir-dim}(\widehat{\mc{M}}_{C_p})=k(2-2g)-2\delta(C_p)<0.
\]
If $J$ is generic, then this is impossible by
Lemma~\ref{lem:transversality}(a).
\end{proof}

For the $J$ given by Lemma~\ref{lem:last}, the rest of \S\ref{sec:GFL}
and \S\ref{sec:limits} carries over unchanged, and so the proof of
Theorem~\ref{thm:compactness}(b) is now complete.

\section{The Euler characteristic of flow lines}
\label{sec:euler}

The following theorem is useful for determining what the flow lines
can look like in specific examples.  We use the notation of
\S\ref{sec:begin}.

\begin{theorem}[Euler characteristic]
\label{thm:euler}
Let $C$ be a flow line without trivial cylinders such
that $(\phi,J)$ is admissible near its ends, and let $\tau$ be any
trivialization of $E$ over the ends.  Then:
\begin{description}
\item{(a)}
$\chi(C) \ge c_\tau(C)-\mu_\tau(C)+\mu_\tau^0(C)-Q_\tau(C)$.
\item{(b)}
If $J$ is generic and $\murel(C)=1$, then equality holds in (a).
\end{description}
\end{theorem}

\begin{proof}
(a) The inequality follows by combining the relative
adjunction formulas \eqref{eqn:adjunctionA} and
\eqref{eqn:adjunctionB} and the inequality \eqref{eqn:braids}.

(b) If $J$ is generic, then by the relative adjunction formulas and
Lemma~\ref{lem:indexBeginning},
\[
\chi(C)-c_\tau(C)+\mu_\tau(C)-\mu_\tau^0(C)+Q_\tau(C)\le 
\murel(C)-\dim(\mc{M}_C).
\]
If $\murel(C)=1$, then since $C$ is not a union of trivial cylinders,
$\R$ acts nontrivially on $C$ by translation on $\R\times Y$, so
$\dim(\mc{M}_C)\ge 1$, so equality holds in (a).
\end{proof}

\section{Concluding remarks}
\label{sec:conclusion}

1. The combinatorics in the proof of the index inequality
   \eqref{eqn:indexInequality} are delicate, because the inequality is
   sharp, at least in the following sense: If $C$ is an admissible
   flow line without trivial cylinders, and if the leading coefficient
   in the asymptotic expansion of each end is nonzero, cf.\ equation
   \eqref{eqn:expansion}, then equality holds.  Some key points in
   proving this are as follows.  The condition on the leading
   coefficients implies that equality holds in the winding bound
   \eqref{eqn:windingBound}.  Admissibility then implies that equality
   also holds in the inequalities \eqref{eqn:cabling} and
   \eqref{eqn:linking}, because
   $\op{gcd}(\lceil\mu_\tau(\gamma^q)/2\rceil,q)=1$ for
   $q\in\pin(\gamma,n)$.  Equality then holds in
   Lemma~\ref{lem:incomingInequality}; admissibility implies that
   equality holds in Lemma~\ref{lem:workhorse} where it is applied.
   Likewise equality holds in Lemma~\ref{lem:outgoingInequality}.

2. For the dimensions of moduli spaces of quasi-embeddings that are
not embeddings, an analogue of the index inequality
\eqref{eqn:indexInequality} holds and is strict, by equation
\eqref{eqn:mcg} with $r=d_1=1$, as long as we assume (a) or (b) in
Lemma~\ref{lem:transversality} to ensure transversality.  Roughly
speaking, smoothing a singularity increases the dimension of the
moduli space by a positive even integer.

3. If we replace $\Sigma$ with a symplectic manifold of dimension
greater than 2, then one cannot in general bound the dimension of a
moduli space of embedded curves from above by an additive relative
index defined only in terms of the ends of the flow lines and their
relative homology class.

4. In this paper we studied the mapping torus $Y$, the periodic orbits
of $\phi$, and the symplectic 4-manifold $(\R\times Y,\omega+ds\wedge
dt)$.  One could instead consider a 3-manifold $Y$ with a contact
1-form $\alpha$, the closed orbits of the associated Reeb flow, and
the symplectization $(\R\times Y,d(e^s\alpha))$.  One can consider
embedded pseudoholomorphic curves in the symplectization with ends
asymptotic to Reeb orbits \cite{sft}.  One can state a formal analogue
of our index inequality \eqref{eqn:indexInequality}, and we conjecture
that it is true.  The main issue in proving it, at least following our
approach, would be to better understand the asymptotic behavior of the
ends of flow lines.  Or one might avoid this issue by perturbing the
setup near short Reeb orbits as in \S\ref{sec:nice}.  One should then
further be able to obtain analogues of our compactness results; the
fiber bubbling which can happen for mapping tori has no analogue in
this setting.

\paragraph{Acknowledgments.}
It is a pleasure to thank Y. Eliashberg, D. Salamon, and M. Thaddeus
for many invaluable discussions; the ETH Z\"{u}rich, the Institute for
Theoretical Physics in Santa Barbara, and the Max Planck Institute for
Mathematics in Bonn for their hospitality while parts of this work
were done; and the National Science Foundation for partial financial
support.


\begin{thebibliography}{99}

\bibitem{abbas} C. Abbas, {\em finite energy surfaces and the chord
problem}, Duke Math. J. {\bf 96} (1999), no.\ 2, 241--316.

\bibitem{csm} R. Carter, G. Segal, and I. Macdonald, {\em Lectures on
Lie groups and Lie algebras}, London Math. Soc. Student Texts 32,
Cambridge, 1995.

\bibitem{sft} Ya.\ Eliashberg, A.\ Givental, and H.\ Hofer, {\em
Introduction to symplectic field theory}, math.SG/0010059.

\bibitem{ehwz} Ya.\ Eliashberg, H.\ Hofer, K.\ Wysocki, and E.\
Zehnder, {\em Symplectic field theory I. Compactness}, in preparation.

\bibitem{fhs} A. Floer, H. Hofer, and D. Salamon, {\em Transversality
in elliptic Morse theory for the symplectic action}, Duke
Math. J. {\bf 80} (1995), no. 1, 251--292.

\bibitem{hl} M. Hutchings and Y-J. Lee, {\em Circle-valued Morse
theory, Reidemeister torsion, and Seiberg-Witten invariants of
3-manifolds}, Topology {\bf 38} (1999), no.\ 4, 861--888.

\bibitem{ht1} M. Hutchings and M. Thaddeus, {\em Periodic Floer
homology}, in preparation.

\bibitem{hwz1} H. Hofer, K. Wysocki, and E. Zehnder, {\em Properties
of pseudoholomorphic curves in symplectisations. I. Asymptotics},
Ann. Inst. H. Poincar\'{e} Anal. Non Lin\'{e}aire {\bf 13} (1996),
no. 3, 337--379.  (Correction in Ann. Inst. H. Poincar\'{e}
Anal. Non Lin\'{e}aire {\bf 15} (1998), no. 4, 535--538.)

\bibitem{hwz2} H. Hofer, K. Wysocki, and E. Zehnder, {\em Properties
of pseudoholomorphic curves in symplectisations. II. Embedding controls
and algebraic invariants}, Geom. and Func. Anal. {\bf 5}, no. 2
(1995), 270--328.

\bibitem{hwz3} H. Hofer, K. Wysocki, and E. Zehnder, {\em Properties
of pseudoholomorphic curves in symplectizations.  III.  Fredholm
theory}, Topics in nonlinear analysis, 381--475, Progr.\ Nonlinear
Differential Equations Appl., {\bf 35}, Birkh\"{a}user, Basel, 1999.

\bibitem{ip} E. Ionel and T. Parker, {\em Gromov invariants and
symplectic maps}, Math.\ Ann.\ {\bf 314} (1999), no.\ 1, 127--158.

\bibitem{ip2} E. Ionel and T. Parker, {\em The symplectic sum formula
for Gromov-Witten invariants\/}, math.SG/0010217.

\bibitem{k} T. Kato, {\em Perturbation theory for linear operators},
Classics in Mathematics reprint, Springer-Verlag, 1995.

\bibitem{mcduff} D. McDuff, {\em Singularities and positivity of
intersections of $J$-holomorphic curves}, pp. 191--216 in {\em
Holomorphic curves in symplectic geometry} (M. Audin and
F. Lafontaine, ed.), Progress in Mathematics 117, Birkh\"{a}user, 1994.

\bibitem{ms} D. McDuff and D. Salamon, {\em J-holomorphic
curves and quantum cohomology}, University Lecture Series 6, Amer.\
Math.\ Soc., 1994.

\bibitem{sa} D. Salamon, {\em Seiberg-Witten invariants of mapping
tori, symplectic fixed points, and Lefschetz numbers}.  Proceedings of
6th G\"{o}kova Geometry-Topology Conference.  Turkish J. Math. {\bf
23} (1999), no. 1, 117--143.

\bibitem{sz} D. Salamon and E. Zehnder, {\em Morse theory for periodic
solutions of Hamiltonian systems and the Maslov index}, Comm.\ Pure
Appl.\ Math {\bf 45} (1992), no.\ 10, 1303--1360.

\bibitem{s} M. Schwarz, {\em Cohomology operations from
$S^1$-cobordisms in Floer homology}, PhD thesis, ETH Z\"{u}rich, 1995.

\bibitem{t1} C. H. Taubes, {\em The Seiberg-Witten and
Gromov invariants}, Math. Res. Lett. {\bf 2} (1995), no. 2, 221--238.

\bibitem{t2} C. H. Taubes, {\em Counting
pseudo-holomorphic submanifolds in dimension four}, J. Differential
Geom. {\bf 44} (1996), no. 4, 818--893.

\bibitem{t98} C. H. Taubes, {\em The structure of pseudoholomorphic
subvarieties for a degenerate almost complex structure and symplectic
form on $S^1\times B^3$}, Geom. Top. {\bf 2} (1998), 221--332.

\bibitem{t3} C. H. Taubes, {\em The Seiberg-Witten invariants and
4-manifolds with essential tori}, Geom. Topol. {\bf 5} (2001), 441--519.

\bibitem{y} C. Young, {\em A new proof of Gromov compactness in
dimension 4}, in preparation.

\end{thebibliography}
\end{document}